\documentclass[11pt,oneside,reqno]{amsart}

\usepackage[margin=1in]{geometry}
\usepackage{amsmath,amssymb,amsthm,textcomp}
\usepackage{amsfonts,graphicx}
\usepackage{amsaddr}
\usepackage[mathscr]{eucal}
\usepackage{pgfplots}
\pgfplotsset{compat=1.15}
\usepackage{subfig,tikz}
\usepackage{array}
\usepackage{booktabs}
\usepackage{multicol}
\usepackage{multirow}
\usepackage[utf8]{inputenc}
\usepackage{bm}

\usepackage{color}
 
\usepackage{float}
\usepackage{hyperref}
 \hypersetup{colorlinks,citecolor=blue}
\hypersetup{colorlinks=true,linkcolor=red,filecolor=magenta,    urlcolor=blue}

\usepackage{pslatex}
\usepackage{array}
\usepackage[square,numbers]{natbib}
\allowdisplaybreaks[1]
\theoremstyle{definition}

\newcommand{\ncom}{\newcommand}

\ncom{\integ}[4]{\int_{#1}^{#2}\,{#3}\,d{#4}}
\ncom{\vspan}[1]{{{\rm\,span}\{ #1 \}}}
\ncom{\dm}[1]{ {\displaystyle{#1} } }
\ncom{\ri}[1]{{#1} \index{#1}}

\newtheorem{theorem}{\bf Theorem}[section]
\newtheorem{remark}{\bf Remark}[section]
\newtheorem{proposition}{Proposition}[section]
\newtheorem{lemma}{Lemma}[section]
\newtheorem{corollary}{Corollary}[section]

\newtheorem{definition}{Definition}[section]
\newtheoremstyle
{remarkstyle}
{}
{11pt}
{}
{}
{\bfseries}
{:}
{     }
{\thmname{#1} \thmnumber{#2} }

\theoremstyle{remarkstyle}

\DeclareFontFamily{U}{BOONDOX-calo}{\skewchar\font=45 }
\DeclareFontShape{U}{BOONDOX-calo}{m}{n}{
  <-> s*[1.05] BOONDOX-r-calo}{}
\DeclareFontShape{U}{BOONDOX-calo}{b}{n}{
  <-> s*[1.05] BOONDOX-b-calo}{}
\DeclareMathAlphabet{\mathcalboondox}{U}{BOONDOX-calo}{m}{n}
\SetMathAlphabet{\mathcalboondox}{bold}{U}{BOONDOX-calo}{b}{n}
\DeclareMathAlphabet{\mathbcalboondox}{U}{BOONDOX-calo}{b}{n}

\usepackage{tabularx}

\begin{document}
\color{black}       
\title{Time-changed Compound Poisson Process}



\author{Kartik Tathe$^1$ and Sayan Ghosh$^2$}
\address{Department of Mathematics, Birla Institute of Technology and Science, Pilani, Hyderabad Campus, Jawahar Nagar, Kapra Mandal, Medchal District, Telangana 500078, India.}
\thanks{$^2$Corresponding author}
\email{$^1$kartikvtathe@gmail.com,$^2$sayan@hyderabad.bits-pilani.ac.in}





\begin{abstract}
In this paper, we introduce a time-changed compound Poisson process (TCCPP) obtained by subordinating the classical compound Poisson process (CPP) with a subordinator and an inverse subordinator, which are defined via distinct Bernstein functions. We study the governing equations and distributional properties of a TCCPP including its moments, dependence structure and limiting behaviour. We derive distributions of the integrated CPP and a linear combination of TCCPPs. We also obtain the law of iterated logarithm, the central limit theorem and weighted sum representation for a TCCPP. Finally, we provide some characterizations of a TCCPP in terms of martingales, a univariate multinomial process, a triangular array process and waiting time distribution. The proposed framework generalizes several existing fractional stochastic processes and offers greater flexibility for modelling complex count data.
\end{abstract}

\maketitle


\section{Introduction} 
The Poisson process plays a central role in stochastic modeling due to its simplicity and ability to describe random counting phenomena with independent and stationary increments. However, in many real world situations, events do not occur uniformly over time. Instead, they may occur in clusters with varying rates and also show dependence between different time periods. In addition, although the Poisson distribution is fundamental, its equidispersion property, with variance equal to mean, often fails in practice. As a result, the Poisson process is often inadequate to describe real world phenomena. A natural and widely studied generalization that addresses some of these limitations is the compound Poisson process (CPP), in which each arrival is associated with a random batch size. This leads to a more flexible model for aggregate quantities and has found wide applications in areas such as insurance risk theory, finance, queueing systems and reliability theory.

Despite its flexibility, the classical CPP still inherits key limitations from the underlying Poisson process. In particular, it assumes exponential inter-arrival times and hence exhibits the memoryless property, that is, it satisfies the Markov property. In recent years, many efforts have been made to overcome these limitations through different generalizations. Among them, fractional or time-changed approaches have been useful in capturing long or short range dependence, heavy-tailed waiting times and irregular behavior. This is typically achieved by introducing a random time change through an independent stochastic clock leading to processes with non-Markovian behavior. 

The first fractional extension of the CPP was introduced by \citet{Mainardi2004} and \citet{Meerschaert2004}. A semi-Markov version was later examined by \citet{Scalas2011}. Subsequently, \citet{Beghin2009} and \citet{Beghin2012} developed alternative formulations of fractional compound Poisson processes and studied their main properties, including time-changed representations, governing equations and limiting behavior. These models are typically constructed by replacing the time parameter $t$ with an independent inverse stable subordinator and, in some cases, by modifying the jump size distribution (see \citet{Beghin2014}). Later, a multivariate version of such fractional models was proposed by \citet{Beghin2016}. In \citet{Crescenzo2015}, a CPP subordinated by a Poisson process was analyzed. More recently, some specific fractional variants of the CPP have been explored in \citet{Singh2020}, \citet{Gupta2023}, \citet{Khandakar2023}, \citet{KatariaGFCP}, and \citet{Gupta2026}.

However, in many practical situations, jumps or arrivals do not occur at constant rates, and their magnitudes or batch sizes may also vary over time. Moreover, the jumps or arrivals at an instant may depend on the same in some previous instances. Such features often arise in queueing systems and traffic flow. These observations motivate the need for more flexible models compared to a CPP. In this article, we study a time-changed compound Poisson process (TCCPP) driven by a subordinator with Bernstein function $f$ and an inverse subordinator with Bernstein function $g$. No specific forms of the subordinators and jump size distribution of the proposed process are assumed. This allows us to develop a more unified and general framework compared to the existing processes in the literature.   

The rest of the paper is organized as follows. Section \ref{sec 2} presents the preliminaries and definitions required for the study. Section \ref{sec 3} introduces the TCCPP. We derive the governing equations for the state probabilities and investigate its distributional properties including mean, variance, covariance and correlation or dependence structure. In Section \ref{sec-4}, we obtain the distribution of a linear combination of TCCPPs along with the fractional integral of a TCCPP, weighted sum representation of a TCCPP and order statistics related to a TCCPP. Section \ref{sec-5} explores the asymptotic behaviour of a TCCPP and its related processes including the law of iterated logarithm and a central limit theorem for the TCCPP. Finally, in Section \ref{sec-6}, we provide a martingale characterization of the TCCPP and also characterizations based on a univariate multinomial process, triangular array process and waiting time distribution. 
\section{Preliminaries} \label{sec 2}
 In this section, we present some preliminary results and definitions which will be used later in the paper.
\subsection{Compound Poisson process}\label{CPP}
A compound Poisson process (CPP) denoted by $\{\mathbb{M}(t)\}_{t\ge 0}$ is defined as
\begin{equation*}
    \mathbb{M}(t)=\sum_{i=1}^{N(t)}X_i,
\end{equation*}
where $X_i\,,\,\,i=1,2,\ldots$ are i.i.d. jumps with common density $f_X$ and are independent of the Poisson process $\{N(t)\}_{t\ge 0}$ with intensity $\lambda>0$. The  p.m.f. of the CPP is given by
\begin{equation}\label{pmf_CPP}
    \mathbb{P}(\mathbb{M}(t)=n)=p(n,t)=\sum_{m=0}^{\infty}f_X^{*m}(n)\frac{e^{-\lambda t}(\lambda t)^{m}}{m!},
\end{equation}
where $f_X^{*m}$ is the $m$-fold convolution of the density $f_X$. 
\subsection{Subordinator}\label{subordinator}
A subordinator $\{D_f(t)\}_{t\geq 0}$ is a one dimensional L\'evy process whose Laplace transform (see \citet{Applebaum2009}, Section 1.3.2) is given by
\[
\mathbb{E}[e^{-sD_f(t)}]
=e^{-t f(s)}, \quad s\geq 0,
\]
where 
\[
f(s)=bs+\int_0^\infty (1-e^{-sx})\,\nu(dx),
\]
is called the Bernstein function or Laplace exponent. Here $b\geq 0$ is the drift coefficient and $\nu$ is a non-negative L\'evy measure on $[0,\infty)$ satisfying $\nu([0,\infty))=\infty$ and $\int_0^\infty (1\wedge x)\,\nu(dx)<\infty$. The sample paths of a subordinator $\{D_f(t)\}_{t\geq 0}$ are non-decreasing with $D_f(0)=0$ a.s.
\subsection{Inverse subordinator}\label{inverse_subordinator}
Let $\{Z_g(t)\}_{t\geq 0}$ be a subordinator with Bernstein
function $g$. The inverse subordinator $\{Y_g(t)\}_{t\geq 0}$
is defined by
\[
Y_g(t)=\inf\{u\geq 0:Z_g(u)>t\}, \quad t\geq 0.
\]
It represents the first passage time of the subordinator
above the level $t$.
\subsection{Generalized R-L and C-D derivatives}\label{RL_CD_relation}
The generalized Caputo-Djrbashian (C-D) derivative for a Bernst\'ein function $f$ is defined as
\begin{equation*}
    \mathcal{D}_t^fu(t)=b\frac{d}{dt}u(t)+\int\frac{\partial}{\partial t}u(t-s)\bar{\nu}(s)ds,
\end{equation*}
where $\bar{\nu}(s)=a+\bar{\nu}(0,\infty)$ is the tail of $\nu$. The generalized Riemann-Liouville (R-L) derivative of $f$ is defined as
\begin{equation*}
    \mathbb{D}_t^fu(t)=b\frac{d}{dt}u(t)+\frac{d}{dt}\int u(t-s)\bar{\nu}(s)ds.
\end{equation*}
The generalized R-L derivative and the generalized C-D derivative satisfy the following relation (see \citet{Toaldo2015})
\begin{equation*}
    \mathbb{D}_t^fu(t)=\mathcal{D}_t^f u(t)+u(0)\bar{\nu}(t).
\end{equation*}
\subsection{LRD and SRD properties}\label{LRD/SRD} For fixed $s>0$, suppose the correlation function of a stochastic process $\{Y(t)\}_{t\geq 0}$ has the following asymptotic behavior (see \citet{DOVIDIO2014} or \citet{Maheshwari2016}): 
\begin{equation*}
\text{Corr}(Y(s), Y(t))\sim c(s)t^{-\theta} \quad\text{as} \quad t\rightarrow \infty \quad \text{for some} \quad c(s)>0.
\end{equation*}
Then $\{Y(t)\}_{t\geq 0}$ has the long range dependence (LRD) property if $\theta\in (0,1)$ and the short range dependence (SRD) property if $\theta\in (1,2)$. 

\section{Time-changed compound Poisson process}\label{sec 3}
Let $\{D_f(t)\}_{t\ge0}$ and $\{Y_g(t)\}_{t\ge0}$ denote a subordinator and an inverse subordinator associated with the Bernstein functions $f$ and $g$ respectively. We assume the subordinators are independent of each other and their moments are finite. It's clear that $\{N_f^g(t)=N(D_f(Y_g(t)))\}_{t\ge 0}$ and $\{N_g^f(t)=N(Y_g(D_f(t)))\}_{t\ge 0}$ represent time-changed Poisson processes where $\{N(t)\}_{t\ge 0}$ is assumed to be independent of $\{D_f(Y_g(t)))\}_{t\ge 0}$ and $\{Y_g(D_f(t)))\}_{t\ge 0}$.  
\begin{definition}\label{def:TCCPP-1}
    We define the time-changed compound Poisson process (TCCPP) as
    \begin{equation*}
        \mathbb{M}_f^g(t):=\sum_{i=1}^{N_f^g(t)}X_i ,
    \end{equation*}
    where $\{X_i\}_{i\ge1}$ are i.i.d. jumps with common distribution $f_X$ and are independent of $\{N_f^g(t)\}_{t\ge 0}$. 
\end{definition}
Note that $\{D_f(t)\}_{t\ge 0}$ has independent and stationary increments. However, the inverse subordinator $\{Y_g(t)\}_{t\ge 0}$ has dependent and non-stationary increments (see Theorem 3.1 and Corollary 3.3 of  \citet{Meerschaert2004}). Hence the process $\{D_f(Y_g(t))\}_{t\ge 0}$ also lacks independent or stationary increments and consequently $\{\mathbb{M}_f^g(t)=\mathbb{M}(D_f(Y_g(t)))\}_{t\ge 0}$ is not a L\'evy process. Alternatively, observe that $\{\mathbb{M}_f^g(t)\}_{t\ge 0}$ is not a L\'evy process by definition since $\{N_f^g(t)\}_{t\ge 0}$ is not a L\'evy process.
\begin{definition}\label{def:TCCPP-2}
    Interchanging $D_f(t)$ and $Y_g(t)$ in Definition \ref{def:TCCPP-1}, we can define TCCPP as
    \begin{equation*}
        \mathbb{M}_g^{f}(t):=\sum_{i=1}^{N_g^f(t)}X_i ,
    \end{equation*}
    where $X_i$'s are i.i.d. jumps independent of $\{N_g^f(t)\}_{t\ge 0}$ with common distribution $f_X$.
\end{definition}
Note that $\{\mathbb{M}_g^{f}(t)=\mathbb{M}(Y_g(D_f(t)))\}_{t\ge 0}$ is not a L\'evy process using similar arguments as above. Henceforth, we shall consider Definition \ref{def:TCCPP-1} of a TCCPP unless otherwise mentioned. 

\subsection{Special cases of TCCPP}\label{special_cases_TCCPP}
Several counting processes in the literature can be obtained as special cases of the TCCPP by specifying the time change (of the underlying Poisson process) and/or the distribution of the jump sizes given by $f_X$. For example, if $Y_g(t)=t$ a.s., then the TCCPP reduces to the fractional generalization of the CPP studied in \citet{Gupta2026}. Setting $D_f(Y_g(t))=N_2(t)$ a.s., where $\{N_2(t)\}_{t\geq0}$ is a Poisson process independent of $\{N(t)\}_{t\geq0}$, the TCCPP reduces to a compound Poisson process subordinated by a Poisson process, as studied in \citet{Crescenzo2015}. If $N(D_f(Y_g(t)))$ is replaced by $M(t)$ where $\{M(t)\}_{t\ge 0}$ is a generalized counting process (GCP) (see \citet{Crescenzo2016}), then the TCCPP is called the compound GCP (see \citet{Khandakar2024}). Moreover if $X_i\overset{d}{=} M(1)$ and $N(D_f(Y_g(t)))$ is replaced by $M_0(t)$ where $\{M_0(t)\}_{t\ge 0}$ is independent of $\{M(t)\}_{t\ge 0}$, then the TCCPP has the same distribution as an iterated GCP (see \citet{Dhillon2026}). For all special cases of the TCCPP listed below, we assume $f_X=\mathbb{P}(X_i=j)=\frac{\lambda_j}{\Lambda}$ for $1\le j\le k$ ($\Lambda=\sum_{j=1}^k\lambda_j$), that is, a jump of size $j$ of the TCCPP occurs with rate or intensity $\lambda_j$. 

First, if $D_f(Y_g(t))=t$ a.s., then the TCCPP reduces to the GCP. Further, if $\lambda_1=\cdots =\lambda_k$, then the GCP reduces to the Poisson process of order $k$ (PPoK) (see \citet{Kostadinova2013}) and to the compound Poisson process of order $k$ (CPPoK) (see \citet{Singh2020}) if $\{N(t)\}_{t\ge 0}$ is replaced by PPoK. Note that PPoK reduces to the classical Poisson process for $k=1$. Next, if $Y_g(t)=t$ a.s. and $\lambda_1=\cdots=\lambda_k$, the TCCPP reduces to the generalized space-fractional Poisson process of order $k$ (GSFPPoK) (see \citet{Gupta2023})

Taking $D_f(t)=D_\beta(t)$ ($\beta$-stable subordinator) a.s. and $Y_g(t)=Y_\alpha(t)$ (inverse $\alpha$-stable subordinator) a.s., the TCCPP reduces to the generalized space-time fractional counting process (GSTFCP) (see \citet{Tathe2025}). For $\beta=1$ and $\alpha=1$, GSTFCP reduces to the generalized space-fractional counting process (GSFCP) and the generalized fractional counting process (GFCP) (see \citet{KatariaGFCP}), respectively. Moreover, if $\lambda_1=\cdots=\lambda_k$, GSTFCP reduces to the space-time fractional Poisson process of order $k$ (STFPPoK), which in turn reduces to FPPoK for $\beta=1$ and to SFPPoK for $\alpha=1$. Observe that GSTFCP, GSFCP and GFCP reduce to STFPP, SFPP (see \citet{Orsingher2012}) and TFPP (see \citet{LASKIN2003}) respectively for $k=1$.

Now if $Y_g(t)= Y_{\beta}(t)$ a.s., then the process $\{\mathbb{M}_f^g(t)\}_{t\ge 0}$ in Definition \ref{def:TCCPP-2} reduces to TCGFCP-I (see \citet{Kataria2022b}), which in turn reduces to TCGCP-I for $\beta=1$. On the other hand, if $Y_g(t)= Y_{\alpha}(t)$ a.s, then the process $\{\mathbb{M}_g^f(t)\}_{t\ge 0}$ in Definition \ref{def:TCCPP-1} reduces to the processes considered in Sections 6 and 5 (for $\alpha=1$) of \citet{Katariaarxiv}. 
Finally, if $D_f(Y_g(t))=t$ a.s. and the jump sizes are natural numbers (may be infinite) such that the sum of their intensities is finite, then the TCCPP is equal in distribution to a Generalized Poisson Process (see \citet{Vellaisamy2026}).

\subsection{Governing Equations}
\begin{theorem}
    The pmf $p_f^g(n,t)=\mathbb{P}(\mathbb{M}_f^g(t)=n)$ satisfies the following fractional differential equation
    \begin{equation*}
        \mathcal{D}_t^fp_{f}^{g}(n,t)=-\lambda p_{f}^{g}(n,t)+\lambda\int_{-\infty}^{\infty}p_{f}^{g}(n-x,t)f_X(x)dx.
    \end{equation*}
\end{theorem}
\begin{proof}
    Let $h_f^g(y,t)$ be the probability density function of the the subordinator $\{D_f(Y_g(t))\}_{t\ge 0}$ and $p(n,t)$ be pmf of CPP. Using conditional argument, we have
    \begin{align*}
        \mathbb{P}(\mathbb{M}_f^g(t)\in dz)&=\int_{0}^{\infty}\mathbb{P}(\mathbb{M}(x)\in dz)h_f^g(x,t)dx.
    \end{align*}
    Now applying the Generalized Riemann-Liouville derivative on both sides of above equation, we get
    \begin{align}
        \mathbb{D}_t^f\mathbb{P}(\mathbb{M}_f^g(t)\in dz)
        &=\int_{0}^{\infty}p(n,x)\mathbb{D}_t^f h_f^g(u,t)dx\nonumber\\
        &=-\int_{0}^{\infty}p(n,x) \frac{\partial}{\partial x}h_f^g(u,t)dx~(\because \mathbb{D}_{t}^{f}h_{f}^{g}(u,t)=-\frac{\partial}{\partial u}h_{f}^{g}(u,t) (\text{see \citet{Meerschaert2013})})\nonumber\\
        &=-p(n,x)h_f^g(x,t)|_{0}^{\infty}+\int_{0}^{\infty}\frac{\partial}{\partial x}p(n,x)h_f^g(x,t)dx.\label{gov_eq_1}
    \end{align}
    It is known that the distribution of CPP (see \citet{Meerschaert2011}) satisfies the following partial differential equation
    \begin{equation}
        \frac{\partial}{\partial x}p(n,x)=-\lambda p(n,x)+\lambda\int_{-\infty}^{\infty}p(n-u,x)f_X(u)du.\label{CPP_gov_eq}
    \end{equation}
    Substituting \eqref{CPP_gov_eq} in \eqref{gov_eq_1}, we have
    \begin{align}
        \mathbb{D}_t^fp_f^g(n,t)&=p(n,0)h_f^g(0,t)+\int_{0}^{\infty}(-\lambda p(n,x)+\lambda\int_{-\infty}^{\infty}p(n-u,x)f_X(u)du)h_f^g(x,t)dx\nonumber\\
        &=p(n,0)h_f^g(0,t)-\lambda\int_{0}^{\infty}p(n,x)h_f^g(x,t)dx+\lambda\int_0^{\infty}\int_{-\infty}^{\infty}p(n-u,x)f_X(u)h_f^g(x,t)du\,dx\nonumber\\
        &=p(n,0)h_f^g(0,t)-\lambda p_f^g(n,t)+\lambda\int_{-\infty}^{\infty}p_f^g(n-x,t)f_X(x)dx.\nonumber
    \end{align}
    Now, using the relation between generalized Reimann-Liouville derivative $\mathbb{D}_t^f$ and generalized Caputo-Djrbashian derivative $\mathcal{D}_t^f$ (see section \ref{RL_CD_relation}), we obtain the result.
\end{proof}

\begin{remark}
    If $\{\mathbb{M}(t)\}_{t\ge 0}$ is GCP, then the state probabilities $p_f^g(n,t)=\mathbb{P}(\mathbb{M}(D_f(Y_g(t)))=n)$ satisfy the following system of differential equations which can also be obtained by time changing Proposition 4 of \citet{Katariaarxiv} with an independent inverse subordinator $\{Y_g(t)\}_{t\ge 0}$:
    \begin{equation*}
        \frac{d}{dt}p_f^g(n,t)=-f(\Lambda)p_f^g(n,t)-\sum_{m=0}^{n}\sum_{\Omega(k,m)}f^{(z_k)}(\Lambda)(\prod_{j=1}^{k}\frac{(-\lambda_j)^{x_j}}{x_j!})p_f^g(n-m,t),
    \end{equation*}
   where $z_k=\sum_{j=1}^kx_j$. Using Remark 12 in \citet{Katariaarxiv}, it can be equivalently written as
    \begin{equation}\label{diff_eq}
        \frac{d}{dt}p_f^g(n,t)=-f(\Lambda(I-\frac{1}{\Lambda}\sum_{j=1}^{n\wedge k}\lambda_jB^j))p_f^g(n,t)
    \end{equation}
    where $Bp_f^g(n,t)=p_f^g(n-1,t)$. If we choose $f(s)=s^\beta$, the Bernstein function associated with a $\beta$-stable subordinator, then \eqref{diff_eq} reduces to 
    \begin{equation}
        \frac{d}{dt}p_{f}^g(n,t)=-\Lambda^\beta(I-\frac{1}{\Lambda}\sum_{j=1}^{n\wedge k}\lambda_jB^j)^{\beta}p_{f}^g(n,t).\nonumber     
    \end{equation}
\end{remark}

\subsection{Mean, Variance and Covariance} The mean, variance, and covariance structure of the TCCPP are derived below. The mean of the process is given by
\begin{align}
\mathbb{E}(\mathbb{M}_f^g(t))=\mathbb{E}(\mathbb{E}(\sum_{i=1}^{N_f^g(t)}X_i|N_f^g(t)))
&=\sum_{k=1}^{\infty}\mathbb{E}(\sum_{i=1}^{k}X_i)\mathbb{P}(N_f^g(t)=k)
=\lambda\mathbb{E}(X_1)\mathbb{E}(D_f(Y_g(t))).\label{mean_TCCPP}
\end{align}
To obtain the variance, we first evaluate the second moment similar to $\eqref{mean_TCCPP}$
\begin{align}
\mathbb{E}((\mathbb{M}_f^g(t))^2)
&=\sum_{k=1}^{\infty}\mathbb{E}((\sum_{i=1}^{k}X_i)^2)\mathbb{P}(N_f^g(t)=k)\nonumber\\
&=\sum_{k=1}^{\infty}\mathbb{E}(\sum_{i=1}^{k}X_i^2)\mathbb{P}(N_f^g(t)=k)+\sum_{k=1}^{\infty}\mathbb{E}(\sum_{i=1}^{k}\sum_{j=1,\,i\ne j}^{k}X_iX_j)\mathbb{P}(N_f^g(t)=k)\nonumber\\
&=\mathbb{E}(X_1^2)\mathbb{E}(N_f^g(t))+(\mathbb{E}(X_1))^2\mathbb{E}(N_f^g(t)(N_f^g(t)-1)).\label{exp_x2}    
\end{align}
Using \eqref{mean_TCCPP} and \eqref{exp_x2} the variance of $\{\mathbb{M}_f^g(t)\}_{t\ge 0}$ is given by
\begin{align}
\mathbb{V}(\mathbb{M}_f^g(t))&=\lambda\mathbb{E}(X_1^2)\mathbb{E}(D_f(Y_g(t))))+\lambda^2(\mathbb{E}(X_1))^2\mathbb{V}(D_f(Y_g(t)))).\label{variance_TCCPP}
\end{align}
To derive the covariance, we assume $s<t$. Then 
\begin{align}
\mathbb{E}(\mathbb{M}_f^g(s)\times \mathbb{M}_f^g(t))
=\mathbb{E}((\sum_{i=1}^{N_f^g(s)}X_i)^2)+\mathbb{E}(\sum_{i=1}^{N_f^g(s)}X_i\times\sum_{i=N_f^g(s)+1}^{N_f^g(t)}X_i).\label{prod_exp1}
\end{align}
Now
\begin{align}
\mathbb{E}(\sum_{i=1}^{N_f^g(s)}X_i\times\sum_{i=N_f^g(s)+1}^{N_f^g(t)}X_i)
&=\sum_{k=1}^{\infty}\sum_{l=1}^{\infty}\mathbb{E}(\sum_{i=1}^{k}X_i\times\sum_{i=k+1}^{l}X_i)\mathbb{P}(N_f^g(s)=k, N_f^g(t)=l)\nonumber\\
&=\sum_{k=1}^{\infty}\sum_{l=1}^{\infty}k(l-k)(\mathbb{E}(X_1))^2\mathbb{P}(N_f^g(s)=k, N_f^g(t)=l)\nonumber\\
&=(\mathbb{E}(X_1))^2\mathbb{E}\{N_f^g(s)(N_f^g(t)-N_f^g(s))\}.\label{prod_exp2}
\end{align}
Hence, combining \eqref{mean_TCCPP}, \eqref{exp_x2}, \eqref{prod_exp1}, and \eqref{prod_exp2}, the covariance of $\{\mathbb{M}_f^g(t)\}_{t\ge 0}$ is given by 
\begin{align}
\mathrm{Cov}(\mathbb{M}_f^g(s),\mathbb{M}_f^g(t))
&=\mathbb{E}(X_1^2)\mathbb{E}(N_f^g(s))+(\mathbb{E}(X_1))^2\mathbb{E}(N_f^g(s)(N_f^g(s)-1))\nonumber\\
    &+(\mathbb{E}(X_1))^2\mathbb{E}\{N_f^g(s)(N_f^g(t)-N_f^g(s))\}-(\mathbb{E}(X_1))^2\mathbb{E}(N_f^g(s))\mathbb{E}(N_f^g(t))\nonumber\\
    &=\mathbb{V}(X_1)\mathbb{E}(N_f^g(s))+(\mathbb{E}(X_1))^2 \mathrm{Cov}(N_f^g(s),N_f^g(t)).\label{cov2}
\end{align}
Using Lemma 2.1 of \citet{Leonenko2014}, we have
\begin{align}
\mathrm{Cov}(N_f^g(s),N_f^g(t))&=\mathbb{V}(N(1))\mathbb{E}(D_f(Y_g(s)))+(\mathbb{E}(N(1)))^2\mathrm{Cov}(D_f(Y_g(s)),D_f(Y_g(t)))\nonumber\\
&=\lambda\mathbb{E}(D_f(Y_g(s)))+\lambda^2\mathrm{Cov}(D_f(Y_g(s)),D_f(Y_g(t))).\label{cov3}
\end{align}
Substituting \eqref{cov3} in \eqref{cov2}, we obtain
\begin{align}
\mathrm{Cov}(\mathbb{M}_f^g(s),\mathbb{M}_f^g(t))&=\lambda\mathbb{E}(X_1^2)\mathbb{E}(D_f(Y_g(s)))+\lambda^2(\mathbb{E}(X_1))^2\mathrm{Cov}(D_f(Y_g(s)),D_f(Y_g(t))).\label{Cov_TCCPP}
\end{align}

\begin{remark}For the TCCPP $\{\mathbb{M}_f^g(t)\}_{t\ge0}$, the index of dispersion is given by
\begin{align}
\mathbb{V}(\mathbb{M}_f^g(t))-\mathbb{E}(\mathbb{M}_f^g(t))&=\lambda\mathbb{E}(X_1^2)\mathbb{E}(D_f(Y_g(t))))+\lambda^2(\mathbb{E}(X_1))^2\mathbb{V}(D_f(Y_g(t))))-\lambda\mathbb{E}(X_1)\mathbb{E}(D_f(Y_g(t)))\nonumber\\
&=\lambda\mathbb{E}(D_f(Y_g(t))))(\mathbb{E}(X_1^2)-\mathbb{E}(X_1))+\lambda^2(\mathbb{E}(X_1))^2\mathbb{V}(D_f(Y_g(t)))).\nonumber
\end{align}
The TCCPP exhibits over-dispersion whenever $\mathbb{V}(\mathbb{M}_f^g(t))>\mathbb{E}(\mathbb{M}_f^g(t))$, that is,  $\displaystyle\frac{\mathbb{E}\!(D_f(Y_g(t)))}{\mathbb{V}\!(D_f(Y_g(t)))}\le\frac{\lambda(\mathbb{E}(X_1))^2}{\mathbb{E}(X_1)-\mathbb{E}(X_1^2)}$
and under-dispersion otherwise. 
\end{remark}
We next examine the correlation structure of a TCCPP. For fixed $s>0$ and $t\to\infty$, assume that $\mathbb{E}(D_f(Y_g(t)))\\ \sim c_1 t^{\alpha},~
\mathbb{V}(D_f(Y_g(t))) \sim c_2 t^{2\alpha},~
\mathrm{Cov}(D_f(Y_g(s)),D_f(Y_g(t))) \sim mc_3(s)$, where $c_1,c_2,m>0$, $c_3(s)>0$ and $\alpha>0$. Using $\eqref{variance_TCCPP}$ and $\eqref{Cov_TCCPP}$, the correlation function of the TCCPP is given by
\begin{align}\label{correlation_TCCPP}
\mathrm{Corr}(\mathbb{M}_f^g(s),\mathbb{M}_f^g(t))
&=\frac{
\lambda \mathbb{E}(X_1^2)\mathbb{E}(D_f(Y_g(s)))
+\lambda^2 (\mathbb{E}X_1)^2 \mathrm{Cov}(D_f(Y_g(s)),D_f(Y_g(t)))
}{
\sqrt{\mathbb{V}(\mathbb{M}_f^g(s))}
\sqrt{
\lambda \mathbb{E}(X_1^2)\mathbb{E}(D_f(Y_g(t)))
+\lambda^2 (\mathbb{E}X_1)^2 \mathbb{V}(D_f(Y_g(t)))
}
}\nonumber\\
&= \frac{
\lambda \mathbb{E}(X_1^2)c_1s^{\alpha}
+\lambda^2 (\mathbb{E}X_1)^2 mc_3(s)}{
\sqrt{\mathbb{V}(\mathbb{M}_f^g(s))}
\sqrt{
\lambda \mathbb{E}(X_1^2)c_1 t^{\alpha}
+\lambda^2 (\mathbb{E}X_1)^2 c_2 t^{2\alpha}
}
} \sim c(s)t^{-\alpha}~\text{as}~t\to\infty,
\end{align}
where 
\begin{equation*}
c(s) = \frac{\lambda \mathbb{E}(X_1^2)c_1s^{\alpha}
+\lambda^2 (\mathbb{E}X_1)^2 mc_3(s)}{\sqrt{\mathbb{V}(\mathbb{M}_f^g(s))}
\sqrt{\lambda^2 (\mathbb{E}X_1)^2c_2}}.
\end{equation*}

From Section \ref{LRD/SRD}, the TCCPP exhibits long-range dependence (LRD) property for $0<\alpha<1$ and short-range dependence (SRD) property for $1<\alpha<2$.
\begin{remark}
Under the assumptions of Section \ref{special_cases_TCCPP}, if the inverse subordinator $\{Y_g(t)\}_{t\ge0}$ is chosen as the inverse $\alpha$-stable subordinator $\{Y_\alpha(t)\}_{t\ge0}$ ($0<\alpha<1$) corresponding to the Bernstein function $g(s)=s^\alpha$, then the resulting process coincides with the time-changed TCGCP-I (see Section 6 of \citet{Katariaarxiv}). Consequently, the dependence structure obtained through \eqref{correlation_TCCPP} extends the dependence structure of the time-changed TCGCP-I (see Theorem 6 of \citet{Katariaarxiv}).
\end{remark}

\section{Properties}\label{sec-4}
In this section, we study some distributional properties related to  a linear combination, the fractional integral, weighted sum representation and order statistics for a TCCPP.
\subsection{Linear combination of TCCPPs}
Consider $k$ TCCPPs denoted by $\displaystyle \mathbb{M}_{f_m}^{g_m}(t)=\sum_{i=1}^{N_{mf_m}^{g_m}(t)} X_i^{(m)},\, 1 \le m \le k$, where for each $m$, $\{X_i^{(m)}\}_{i\ge1}$ is a sequence of non-negative i.i.d.\ random variables with common density function $f_X^{m}$ and also independent of $\{N_{mf_m}^{g_m}(t)))\}_{t\ge 0}$. Assume that the random variables $X_i^{(1)}, X_i^{(2)}, \ldots, X_i^{(k)}$ are independent for each fixed $i$ and the sequences $\{X_i^{(1)}\}_{i\ge1}, \{X_i^{(2)}\}_{i\ge1}, \dots, \{X_i^{(k)}\}_{i\ge1}$ are jointly independent of $\{N_{1f_1}^{g_1}(t), \dots, N_{kf_k}^{g_k}(t)\}_{t\ge 0}$. Define a linear combination of the above TCCCPs as $\displaystyle \mathcal{M}_f^g(t) = \sum_{m=1}^k a_m \mathbb{M}_{f_m}^{g_m}(t)$ where $a_1,\dots,a_k \in \mathbb{R}$.
\begin{theorem}\label{linear_combination_TCCPP}
The moment generating function of $\{\mathcal{M}_f^g(t)\}_{t\ge 0}$ is given by
\[
\mathbb{E}[e^{s \mathcal{M}_f^g(t)}]
=\mathbb{E}(
\exp(\sum_{m=1}^{k}
\lambda_m D_{f_m}(Y_{g_m}(t))
(\mathbb{E}(e^{s a_m X_1^{(m)}})-1)
)
)
\]
and its probability mass function is given by
\begin{align}
\mathbb{P}(\mathcal{M}_f^g(t)=n)
=\int_{0}^{\infty}\cdots\int_{0}^{\infty}
\sum_{\substack{n_1+\cdots+n_k=n}}
\prod_{m=1}^k 
(
\sum_{r_m=0}^{\infty}
f_{X_m}^{*r_m}\!(n_m/a_m)
\frac{e^{-\lambda_m x_m}(\lambda_m x_m)^{r_m}}{r_m!}
)\,
\prod_{m=1}^k h_{f_m}^{g_m}(x_m,t)\,dx_m,\notag
\end{align}
where $f_{X_m}^{*r_m}(n_m/a_m)=0$ whenever $n_m/a_m \notin \mathbb{N}_0$.
\end{theorem}
\begin{proof}
    Using the independence assumptions, we obtain for any $n_1,\dots,n_k\ge 0$,
\begin{align*}
&\mathbb{P}(\mathbb{M}_{f_1}^{g_1}(t)=a_1,\dots,\mathbb{M}_{f_k}^{g_k}(t)=a_k) \\
&=\sum_{n_1=0}^{\infty}\cdots \sum_{n_k=0}^{\infty} 
\mathbb{P}(\sum_{i=1}^{n_1} X_i^{(1)}=a_1,\dots,\sum_{i=1}^{n_k} X_i^{(k)}=a_k)\mathbb{P}(N_{1f_1}^{g_1}(t)=n_1,\dots,N_{kf_k}^{g_k}(t)=n_k)\\
&=\sum_{n_1=0}^{\infty}\cdots \sum_{n_k=0}^{\infty} 
\mathbb{P}(\sum_{i=1}^{n_1} X_i^{(1)}=a_1)\dots\mathbb{P}(\sum_{i=1}^{n_k} X_i^{(k)}=a_k)\mathbb{P}(N_{1f_1}^{g_1}(t)=n_1)\dots\mathbb{P}(N_{kf_k}^{g_k}(t)=n_k)\\
&=\sum_{n_1=0}^{\infty}\mathbb{P}(\mathbb{M}_{f_1}^{g_1}(t)=a_1|N_{1f_1}^{g_1}(t)=n_1)\mathbb{P}(N_{1f_1}^{g_1}(t)=n_1)\times\cdots\times \sum_{n_k=0}^{\infty}\mathbb{P}(\mathbb{M}_{f_k}^{g_k}(t)=a_k|N_{kf_k}^{g_k}=n_k)\mathbb{P}(N_{kf_k}^{g_k}=n_k)\\
&=\mathbb{P}(\mathbb{M}_{f_1}^{g_1}(t)=a_1)\dots\mathbb{P}(\mathbb{M}_{f_k}^{g_k}(t)=a_k)
\end{align*}
which proves that $\{ \mathbb{M}_{f_1}^{g_1}(t)\}_{t\ge 0},\dots, \{\mathbb{M}_{f_k}^{g_k}(t)\}_{t\ge 0}$ are independent of each other. Hence, using \eqref{m.g.f.TCCPP}
the m.g.f. of $\{\mathcal{M}_f^g(t)\}_{t\ge 0}$ is given by
\begin{align}
\mathbb{E}(\exp(s \sum_{m=1}^k a_m \mathbb{M}_{f_m}^{g_m}(t)))=\prod_{m=1}^k
\mathbb{E}(\exp(sa_m \mathbb{M}_{f_m}^{g_m}(t)))&=\prod_{m=1}^k 
\mathbb{E}(
\exp(
\lambda_m D_{f_m}(Y_{g_m}(t))
(\mathbb{E}(e^{s a_m X_1^{(m)}})-1)
)
)\nonumber\\
&=\mathbb{E}(
\exp(\sum_{m=1}^{k}
\lambda_m D_{f_m}(Y_{g_m}(t))
(\mathbb{E}(e^{s a_m X_1^{(m)}})-1)
)
).\nonumber
\end{align}
To derive the probability mass function of $\{\mathcal{M}_f^g(t)\}_{t\ge 0}$, we first consider the distribution of each component in the sum. Note that the probability mass function of the TCCPP $\{\mathbb{M}_f^g(t)\}_{t\ge 0}$ is 
\begin{align}
p_f^g(n,t)=\mathbb{P}(\mathbb{M}_f^g(t)=n)&=\int_{0}^{\infty}\mathbb{P}(\mathbb{M}(x)=n)h_f^g(x,t)dx =\int_{0}^{\infty}\sum_{m=0}^{\infty}f_X^{*m}(n)\frac{e^{-\lambda x}(\lambda x)^{m}}{m!} h_f^g(x,t)dx,\nonumber
\end{align}
where $\{\mathbb{M}(t)\}_{t\ge 0}$ is a CPP and $f_X^{*m}$ is the $m$-fold convolution of $f_X$. It follows that
\begin{align}\label{pmf_constant_multiple_TCCPP}
    \mathbb{P}(a_m \mathbb{M}_{f_m}^{g_m}(t)=n)&=\int_{0}^{\infty}\sum_{r_m=0}^{\infty}f_{X_m}^{*r_m}(n/a_m)\frac{e^{-\lambda_m x_m}(\lambda_m x_m)^{r_m}}{r_m!} h_{f_m}^{g_m}(x_m,t)dx_m,~\text{if}~n/a_m\in \mathbb{N}_0.
\end{align}
Since the processes $\{a_m \mathbb{M}_{f_m}^{g_m}(t)\}_{t\ge 0}$ are independent, the distribution of their sum is given by
\begin{align}
\mathbb{P}(\mathcal{M}_f^g(t)=n)=\mathbb{P}(\sum_{m=1}^k a_m \mathbb{M}_{f_m}^{g_m}(t)=n)=\sum_{\substack{n_1+\cdots+n_k=n}}
\prod_{m=1}^k 
\mathbb{P}(a_m \mathbb{M}_{f_m}^{g_m}(x_m)=n_m).\label{pmf_linear_combination}
\end{align}
Substituting \eqref{pmf_constant_multiple_TCCPP} in \eqref{pmf_linear_combination}, we obtain the required p.m.f. of $\{\mathcal{M}_f^g(t)\}_{t\ge 0}$.
\end{proof}
Theorem \ref{linear_combination_TCCPP} shows that a linear combination of independent TCCPPs is not a TCCPP.

\subsection{Fractional integral of the TCCPP}
We consider the Riemann-Liouville fractional integral (see \citet{Orsingher2013}) of the TCCPP $\{\mathbb{M}_f^g(t)\}_{t\ge0}$ defined by
\begin{equation}\label{frac_int_TCCPP}
\mathcal{I}_{\mathbb{M}_f^g}^{\zeta}(t):=\frac{1}{\Gamma(\zeta)}\int_{0}^{t}(t-s)^{\zeta-1}\mathbb{M}_f^g(s)\,ds,\quad t\ge0,\;\zeta>0.\nonumber
\end{equation}
The motivation for studying the fractional integral of the TCCPP stems from the fact that integrated non-negative stochastic processes and more generally integrated counting processes, frequently arise in the applied mathematics literature (see \citet{Jerwood1970}, \citet{Downton1972}, \citet{Hernandez1999}, \citet{Pollett2003}).
For integer values of $\zeta=m\in\mathbb{N}$, the fractional integral admits the iterated integral representation given by
\begin{equation}
\frac{1}{(m-1)!}\int_{0}^{t}(t-s)^{m-1}\mathbb{M}_f^g(s)\,ds=\int_{0}^{t}ds_1\int_{s_1}^{t}ds_2\cdots\int_{s_{m-1}}^{t}\mathbb{M}_f^g(s_m)\,ds_m.\nonumber
\end{equation}
The mean of the fractional integral is 
\begin{align}
\mathbb{E}\!(\mathcal{I}_{\mathbb{M}_f^g}^{\zeta}(t))=\frac{1}{\Gamma(\zeta)}\int_{0}^{t}(t-s)^{\zeta-1}\mathbb{E}\!(\mathbb{M}_f^g(s))ds&=\frac{\lambda\mathbb{E}(X_1)}{\Gamma(\zeta)}
\int_{0}^{t}
(t-s)^{\zeta-1}
E(D_f(Y_g(s)))\,ds \quad(\text{using}\,\eqref{mean_TCCPP})
\nonumber
\end{align}
and for $0\le t_1\le t_2$, the covariance between the fractional integrals at times $t_1$ and $t_2$ is 
\begin{align*}
\mathrm{Cov}(\mathcal{I}_{\mathbb{M}_f^g}^{\zeta}(t_1),\mathcal{I}_{\mathbb{M}_f^g}^{\zeta}(t_2))
&=\frac{1}{(\Gamma(\zeta))^2}\int_0^{t_1}\int_0^{t_2}(t_1-s)^{\zeta-1}(t_2-r)^{\zeta-1}\mathrm{Cov}(\mathbb{M}_f^g(s),\mathbb{M}_f^g(r))\,dr\,ds,
\end{align*}
where $\mathrm{Cov}(\mathbb{M}_f^g(s),\mathbb{M}_f^g(r))$ is given by $\eqref{Cov_TCCPP}$. If $t_1=t_2$, we obtain the variance of the fractional integral.

For $\zeta=1$ and $D_f(Y_g(t))=t$ a.s., $\displaystyle\{\mathcal{I}_{\mathbb{M}_f^g}^{\zeta}(t)\}_{t\ge 0}$ reduces to the integrated compound Poisson process $\displaystyle\{I(t):=\int_0^t \mathbb{M}(s)ds\}_{t\ge 0}$. For fixed $t>0$, the random variable $I(t)$ represents the area under the sample path of $\{\mathbb{M}(t)\}_{t\ge0}$ over the interval $[0,t]$. The integral of the fractional Poisson process was investigated by \citet{Orsingher2013}, who derived a random-sum representation for the integrated Poisson process. The following result establishes an analogous distributional representation for the integrated compound Poisson process, extending the framework from single-event arrivals to clustered arrivals with random jump sizes.
\begin{theorem}
For a CPP $\{\mathbb{M}(t)\}_{t\ge0}$, we have $\displaystyle\int_0^t \mathbb{M}(s)\,ds\stackrel{d}{=}\sum_{k=1}^{N(t)}U_kX_k$, where $\{U_k\}_{k\ge1}$ are i.i.d. $\mathrm{Uniform}[0,t]$ random variables which are jointly independent of $\{X_k\}_{k\ge1}$ and also independent of $\{N(t)\}_{t\ge 0}$.
\end{theorem} 
\begin{proof}
    Let $\{T_n\}_{n\ge1}$ denote the arrival times of the Poisson process $\{N(t)\}_{t\ge0}$. Since $\displaystyle N(t)=\sum_{j\ge 1}\mathbb{I}_{\{T_j\le t\}}$, the CPP $\{\mathbb{M}(t)\}_{t\ge 0}$ can be represented as $\displaystyle \mathbb{M}(t)=\sum_{j=1}^{N(t)}X_j
=\sum_{j\ge1}X_j \mathbb{I}_{\{T_j\le t\}}$. Therefore
\begin{align*}
I(t):=\int_0^t \mathbb{M}(s)ds
=\int_0^t \sum_{j\ge1}X_j \mathbb{I}_{\{T_j\le s\}}ds =\sum_{j\ge1}X_j\int_0^t \mathbb{I}_{\{T_j\le s\}}ds=\sum_{j\ge1}X_j(t-T_j)\mathbb{I}_{\{T_j\le s\}}=\sum_{j=1}^{N(t)}X_j(t-T_j).
\end{align*}
By definition of a CPP, the jump sizes $\{X_j\}_{j\ge1}$ are i.i.d. and independent of the Poisson process $\{N(t)\}_{t\ge0}$. Moreover, the arrival times $\{T_j\}_{j\ge1}$ satisfy $\{T_j\le t\}=\{N(t)\ge j\}$. Hence, each $T_j$ is measurable with respect to $\sigma(\{N(t)\},\,t\ge 0)$. It follows that the sequence $\{X_j\}_{j\ge1}$ is independent of arrival times $\{T_j\}_{j\ge1}$. Consequently, conditional on $(T_1,\ldots,T_{N(t)},N(t))$, the random variables $X_1(t-T_1),\ldots,X_{N(t)}(t-T_{N(t)})$ are independent. Also, $X_j(t-T_j)\perp\!\!\!\!\!\perp N(t)\mid (T_1,\ldots,T_{N(t)})$ for $1\le j\le N(t)$. Using the above facts, the characteristic  function of $\{I(t)\}_{t\ge 0}$ is given by
\begin{align}
\phi_{I(t)}(v)=\mathbb E\!(\exp\!(iv\sum_{j=1}^{N(t)}X_j(t-T_j))
)&=
\mathbb E\!(\mathbb E\!(\exp\!(iv\sum_{j=1}^{N(t)}X_j(t-T_j))\,|\,N(t))) \nonumber\\
&=\mathbb E\!(\mathbb E\!(\mathbb E\!(\exp\!(iv\sum_{j=1}^{N(t)}X_j(t-T_j))\,|\,T_1,\ldots,T_{N(t)},N(t))\,|\,N(t)))\nonumber \\
&=\mathbb E\!(\mathbb E\!(\prod_{j=1}^{N(t)}\mathbb E\!(\exp\!(ivX_j(t-T_j))\,|\,T_1,\ldots,T_{N(t)},N(t))\,|\,N(t))) \nonumber\\
&=\mathbb E\!(\mathbb E\!(\prod_{j=1}^{N(t)}\phi_{X_j}\!(v(t-T_j)\,|\,T_1,\ldots,T_{N(t)})\,|\,N(t)))\nonumber\\
&=\mathbb E\!(\mathbb E\!(\prod_{j=1}^{N(t)}\phi_{X_j}\!(v(t-T_j))\,|\,N(t))).\nonumber
\end{align}
Note that $\displaystyle(T_1,\ldots,T_n)\mid \{N(t)=n\} \stackrel{d}{=}(U_{(1)},\ldots,U_{(n)})$, where $U_{(1)}<\cdots<U_{(n)}$ are the order statistics of $n$ i.i.d. $\mathrm{Uniform}[0,t]$ random variables $U_1,\ldots,U_n$. This implies 
\begin{equation*}
\mathbb E\!(\prod_{j=1}^{N(t)}\phi_{X_j}\!(v(t-T_j))\,|\,N(t))=\mathbb E\!(\prod_{j=1}^{N(t)}\phi_{X_j}\!(v(t-U_{(j)})))=\mathbb E\!(\prod_{j=1}^{N(t)}\phi_{X_j}\!(v(t-U_j)))
\end{equation*}
since $\displaystyle\prod_{j=1}^{n}\phi_{X_j}\!(v(t-u_j))$ is symmetric in $u_1,\ldots,u_n$. 
So the characteristic function of $\{I(t)\}_{t\ge 0}$ takes the form
\begin{equation*}
\phi_{I(t)}(v)=
\mathbb E\!(\mathbb E\!(\prod_{j=1}^{N(t)}\phi_{X_j}\!(v(t-U_j))))=\mathbb E\!(\prod_{j=1}^{N(t)}\phi_{X_j}\!(v(t-U_j)))
=\mathbb E\!(\mathbb E\!(\prod_{j=1}^{N(t)}\phi_{X_j}\!(v(t-U_j)| N(t)))).\nonumber
\end{equation*}
Now
\begin{align}
\mathbb{P}(X_j(t-U_j)\le y\mid N(t))&=\int_0^t\mathbb{P}(X_j(t-U_j)\le y\mid U_j,N(t))\,f_{U_j\mid N(t)}(u_j)\,du_j\nonumber \\
&=\int_0^t\mathbb{P}(X_j(t-U_j)\le y|U_j)\,f_{U_j}(u)\,du_j \quad (\because X_j(t-U_j) \perp\!\!\!\!\!\perp N(t)\mid U_j,~U_j \perp\!\!\!\!\!\perp N(t))\nonumber\\
&=\mathbb{P}(X_j(t-U_j)\le y),\label{independence_X(t-U)_N(t)}
\end{align}
which implies $X_j(t-U_j)$ is independent of $N(t)$. Also 
\begin{equation}
    \mathbb{P}(X_j(t-U_j)\le y)=\frac{1}{t}\int_{0}^{t}\mathbb{P}(X_j\le y/(t-u_j))du_j=\frac{1}{t}\int_{0}^{t}F_{X}(y/(t-u_j))du_j\label{identical_distribution_X(t-U)}
\end{equation} which is independent of $j$, that is, the random variables $X_j(t-U_j)$ are identically distributed. Hence
\begin{align}
\phi_{I(t)}(v)=\mathbb E\!((\mathbb E\!(\phi_{X_1}\!(v(t-U_1))))^{N(t)})&=\mathbb E\!((\frac{1}{t}\int_0^t\phi_{X_1}\!(v(t-s))\,ds)^{N(t)})\nonumber\\
&=e^{-\lambda t}\sum_{n=0}^{\infty}\frac{(\lambda\int_0^t\phi_{X_1}\!(v(t-s))\,ds)^n}{n!} \nonumber\\
&=\exp\!\{\lambda t(\frac{1}{t}\int_0^t\phi_{X_1}(v w)\,dw-1)\}\label{char_fun_LHS},
\end{align}
where the last equality follows from the change of variable $w=t-s$. Note that
\begin{equation}
    \phi_{U_1X_1}(v)=\mathbb{E}(e^{ivU_1X_1})=\mathbb{E}(\mathbb{E}(e^{ivU_1X_1}|U_1))=\mathbb{E}(\phi_{X_1}(vU_1))=\frac{1}{t}\int_{0}^{t}\phi_{X}(vw)dw\label{char_fun_UX}.
\end{equation}
Using \eqref{char_fun_UX} in \eqref{char_fun_LHS}, we have
\begin{equation}
    \phi_{I(t)}(v)=\exp\{\lambda t(\phi_{U_1X_1}(v)-1)\}.\label{char_fun_I(t)}
\end{equation}
Next, replacing $t-U_j$ in \eqref{independence_X(t-U)_N(t)} and \eqref{identical_distribution_X(t-U)} by $U_j$, we infer that $\{U_kX_k\}_{k\ge 1}$ are independent of $\{N(t)\}_{t\ge 0}$ and $\{U_kX_k\}_{k\ge1}$ are identically distributed respectively. Moreover, since $(X_1,\ldots,X_n)$ is independent of $(U_1,\ldots,U_n)$ for $n\in\mathbb{N}$, we have
\begin{align*}
\mathbb{P}(U_1X_1\le y,\ldots,U_nX_n\le y)
&=\frac{1}{t^n}\int_{0}^t\cdots\int_{0}^t\mathbb{P}\!(X_1\le \frac{y}{u_1},\ldots,X_n\le \frac{y}{u_n})du_1\ldots du_n\\
&=\prod_{k=1}^{n}(\frac{1}{t}\int_{0}^t\mathbb{P}\!(X_k\le \frac{y}{u_k})du_k)\\
&=\prod_{k=1}^{n}\mathbb{P}(U_kX_k\le y),
\end{align*}
which shows that $\{U_kX_k\}_{k\ge1}$ are independent random variables. Hence, $\{\sum_{k=1}^{N(t)}U_kX_k\}_{t\ge 0}$ is a CPP with characteristic function $\exp\{\lambda t(\phi_{U_1X_1}(v)-1)\}$. It follows from \eqref{char_fun_I(t)} that 
\begin{equation*}
\int_0^t \mathbb{M}(s)\,ds\stackrel{d}{=}\sum_{k=1}^{N(t)}U_kX_k.
\end{equation*}
\end{proof}

\subsection{Weighted sum representation}
The next result provides a representation of the TCCPP as a weighted sum of independent time-changed Poisson processes where the weights are jump sizes of the TCCPP.
\begin{theorem}\label{Weighted_sum_representation_TCCPP}
The TCCPP with intensity $\lambda$ and jump size distribution $\{p_k\}_{k\ge 1}$ is equal in distribution to the weighted sum of independent time-changed Poisson processes, that is,
\begin{equation*}
    \mathbb{M}_f^g(t) \stackrel{d}{=} \sum_{k=1}^{\infty} k\,N_k(D_f(Y_g(t))),
\end{equation*}
where $\{N_k(D_f(Y_g(t)))\}_{t\ge0}$ are independent time-changed Poisson processes and $N_k(t)$ is independent of $D_{f}(Y_{g}(t))$ for $k\geq 1$ such that $N_k(t)\sim \text{Poi}(\lambda p_k t)$.
\end{theorem}
\begin{proof}
We know that the m.g.f. of a CPP is 
\begin{equation}\label{LT_CPP}
\mathbb{E}(e^{s\mathbb{M}(t)})=\exp(-\lambda t(1-\mathbb{E}(e^{sX_1}))).
\end{equation}
So the m.g.f. of a TCCPP is given by
\begin{align}
\mathbb{E}(e^{s\mathbb{M}_f^g(t)})
=\int_{0}^{\infty}\mathbb{E}(e^{s\sum_{i=1}^{N(x)}X_i})f_{D_f(Y_g(t))}(x)\,dx\notag
&=\int_{0}^{\infty}\exp(\lambda x(\mathbb{E}(e^{sX_1})-1))f_{D_f(Y_g(t))}(x)\,dx \quad (\text{using}~\eqref{LT_CPP})\notag\\
&=\int_{0}^{\infty}\mathbb{E}(e^{\lambda D_f(Y_g(t))(\mathbb{E}(e^{sX_1})-1)}\,|\,D_f(Y_g(t))=x)f_{D_f(Y_g(t))}(x)\,dx\notag\\
&=\mathbb{E}(\mathbb{E}(e^{\lambda D_f(Y_g(t))(\mathbb{E}(e^{sX_1})-1)}\,|\,D_f(Y_g(t))))\notag\\
&=\mathbb{E}(e^{\lambda D_f(Y_g(t))(\mathbb{E}(e^{sX_1})-1)}).\label{m.g.f.TCCPP}
\end{align}

Since the m.g.f. of $X_i$ is $\displaystyle \mathbb{E}(e^{sX_1})=\sum_{k=1}^{\infty} p_k e^{sk}$, we have
$\displaystyle\mathbb{E}(e^{sX_1})-1=\sum_{k=1}^{\infty} p_k (e^{sk}-1)$. Substituting this into \eqref{m.g.f.TCCPP}, we have
\begin{align}
\mathbb{E}(e^{s\mathbb{M}_f^g(t)})&=\mathbb{E}\!(\exp\!(\lambda D_f(Y_g(t)) 
\sum_{k=1}^{\infty} p_k (e^{sk}-1)))=\mathbb{E}\!(\prod_{k=1}^{\infty}
\exp\!(\lambda p_k D_f(Y_g(t))(e^{sk}-1))).\notag
\end{align}
Observe that $\exp\!(\lambda p_k u (e^{sk}-1))$ is the moment generating function of $kN_k(u)$. Now, due to independence of $N_k(D_f(Y_g(t)))$, $k\ge 1$ it follows that
\begin{align}
\mathbb{E}(e^{s\mathbb{M}_f^g(t)})&=\prod_{k=1}^{\infty}
\mathbb{E}\!(e^{skN_k(D_f(Y_g(t)))})=\mathbb{E}(e^{s\sum_{k=1}^{\infty} kN_k(D_f(Y_g(t)))}). \notag
\end{align}
Therefore, by the uniqueness of the moment generating function,
\begin{equation*}
    \mathbb{M}_f^g(t) \stackrel{d}{=} \sum_{k=1}^{\infty} k\,N_k(D_f(Y_g(t))),
\end{equation*}
where $\{N_k(t)\}_{t \ge 0}$ are independent Poisson processes with rates $\lambda p_k$.
\end{proof}

\begin{remark}
Suppose the jump size distribution of a TCCPP is given by $f_X=\mathbb{P}(X_1=j)=\lambda_j/\sum_{j=1}^k\lambda_j$ for $1\le j\le k$. Then the TCCPP reduces to a time-changed GCP denoted by $\{M(D_f(Y_g(t)))\}_{t\ge 0}$. Using Theorem 2.2 of \citet{Dhillon2024}, it follows that $M(D_f(Y_g(t)))=\sum_{j=1}^kjN_j(D_f(Y_g(t)))$ where $\{N_j(t)\}_{t\ge 0}$s are independent Poisson processes with rate $\lambda_j$ and $D_{f}(Y_g(t))$ is independent of $N_j(t)$ for each $j$.    
\end{remark}
\subsection{Order Statistics} Now we investigate the order statistics associated with a random sample whose size is governed by the counting process $\{N_f^g(t)\}_{t\ge0}$.
Let $\{J_i\}_{i\ge 1}$ be a sequence of i.i.d. random variables with common distribution function $F$, independent of the counting process $\{N_f^g(t)\}_{t\ge 0}$. 
Let $\mathbb{I}_{A}$ denote the indicator function of the set $A$. For a fixed $z \in \mathbb{R}$, define the random variable
\begin{equation*}
    B_z(t) := \sum_{i=1}^{N_f^g(t)} \mathbb{I}_{\{J_i < z\}}
\end{equation*}
which counts the number of observations less than $z$ in a sample of size $N_f^g(t)$. Observe that
\begin{equation}\label{binomial counting variable}
    B_{z}(t)|\{N_f^g(t)=n\} ~\sim~\text{Binomial}(n,F(z)).
\end{equation}
For $k \ge 1$, let $J_{(k)}^{f,g}(t)$ denote the $k$-th order statistic of the random sample $\{J_1, J_2, \ldots, J_{N_f^g(t)}\}$. The following result gives the conditional distribution of $J_{(k)}^{f,g}(t)$ given $N_f^g(t) \ge k$.
\begin{theorem}
    Let $\{J_i\}_{i\ge 1}$ be a sequence of i.i.d. random variables with distribution function $F$. Then
    \begin{equation*}
        \mathbb{P}(J_{(k)}^{f,g}(t)\le z|N_f^g(t)\ge k)=\frac{\displaystyle\sum_{n=k}^{\infty}\sum_{j=k}^{n}\binom{n}{j}(F(z))^{j}(1-F(z))^{n-j}\mathbb{P}(N_f^g(t)=n)}{\mathbb{P}(N_f^g(t)\ge k)}.
    \end{equation*}
\end{theorem}\
\begin{proof}
    By definition of conditional probability, 
    \begin{align}
        \mathbb{P}(J_{(k)}^{f,g}(t)\le z|N_f^g(t)\ge k)=\frac{\mathbb{P}(J_{(k)}^{f,g}(t)\le z,~N_f^g(t)\ge k)}{\mathbb{P}(N_f^g(t)\ge k)}
        &=\frac{\displaystyle{\sum_{n=k}^{\infty}}\mathbb{P}(J_{(k)}^{f,g}(t)\le z,~N_f^g(t)=n)}{\mathbb{P}(N_f^g(t)\ge k)}\nonumber\\
        &=\frac{\displaystyle\sum_{n=k}^{\infty}\mathbb{P}(J_{(k)}^{f,g}(t)\le z |~N_f^g(t)=n)\mathbb{P}(N_f^g(t)=n)}{\mathbb{P}(N_f^g(t)\ge k)}.\label{order statistic_eq1}
    \end{align}
    Now
    \begin{align}
        \mathbb{P}(J_{(k)}^{f,g}(t)\le z|~N_f^g(t)=n)&=\mathbb{P}(\text{at least $k$ observations are smaller than $z$ in a sample of size $n$})\nonumber\\
        &=\mathbb{P}(B_z(t)\ge k|~N_f^g(t)= n)\nonumber\\
        &=\sum_{j=k}^{n}\binom{n}{j}(F(z))^{j}(1-F(z))^{n-j}.\quad(\text{using} ~\eqref{binomial counting variable})\label{order statistic_eq2}
    \end{align}
    Using \eqref{order statistic_eq2} in \eqref{order statistic_eq1}, we get
    \begin{align*}
         \mathbb{P}(J_{(k)}^{f,g}(t)\le z|N_f^g(t)\ge k)&=\frac{\displaystyle\sum_{n=k}^{\infty}\sum_{j=k}^{n}\binom{n}{j}(F(z))^{j}(1-F(z))^{n-j}\mathbb{P}(N_f^g(t)=n)}{\mathbb{P}(N_f^g(t)\ge k)}.
    \end{align*}
\end{proof}

\section{Limiting Behaviour}\label{sec-5}
In this section, we investigate the asymptotic behaviour of the TCCPP and related time-changed processes. We first establish distributional limits and a law of the iterated logarithm for the TCCPP. We then derive a limiting distribution for the time-changed Poisson process underlying the TCCPP and establish stable convergence along with a central limit theorem for the TCCPP.
\begin{theorem}
    If $\displaystyle{\lim_{t\rightarrow \infty}\frac{D_f(Y_g(t))}{t^c}\overset{d}{=}Z}$ for some $c>0$, then $\{ \mathbb{M}_f^g(t)\}_{t\ge 0}$ satisfies the following limiting result:
    \begin{equation*}
        \lim_{t\to \infty}\frac{ \mathbb{M}_f^g(t)}{t^c}\overset{d}{=}\lambda \mathbb{E}(X_1)Z.
    \end{equation*}
\end{theorem}
\begin{proof}
    Since $\{X_i\}_{i\ge1}$ are i.i.d.\ random variables with $\mathbb{E}(|X_1|)<\infty$, and $\{N_f^g(t)\}_{t\ge0}$ is an integer-valued random sequence such that $N_f^g(t)\xrightarrow{\mathrm{a.s.}}\infty$ as $t\to\infty$, which is independent of $\{X_i\}_{i\ge1}$, it follows from the Strong Law of Large Numbers for random indices (see \citet{Bertoin2004}, p.~92) that
    \begin{equation}\label{strong_law_eq1}
        \frac{1}{N_f^g(t)}\sum_{i=1}^{N_f^g(t)} X_i \xrightarrow{\mathrm{a.s.}} \mathbb{E}[X_1]\quad \text{as}\quad t\to\infty.
    \end{equation}
    Moreover, since $\{N(t)\}_{t\ge0}$ satisfies $\frac{N(t)}{t}\xrightarrow{\mathrm{a.s.}}\lambda$ as $t\to\infty$, and $D_f(Y_g(t))\xrightarrow{\mathrm{a.s.}}\infty$, it follows that
\begin{equation}\label{strong_law_eq2}
\frac{N_f^g(t)}{D_f(Y_g(t))}\xrightarrow{\mathrm{a.s.}}\lambda \quad \text{as}\quad t\to\infty.
\end{equation}
    Using \eqref{strong_law_eq1} and \eqref{strong_law_eq2}, we obtain
    \begin{align}
        \lim_{t\to \infty}\frac{ \mathbb{M}_f^g(t)}{t^c}&\overset{d}{=}\lim_{t\to \infty}\frac{N_f^g(t)}{(D_f(Y_g(t)))}\frac{(D_f(Y_g(t)))}{t^c}\frac{1}{N_f^g(t)}\sum_{i=1}^{N_f^g(t)} X_i\overset{d}{=}\lambda\mathbb{E}[X_1]Z.\nonumber
    \end{align}
\end{proof}
An analogous result can be established for the process $\{\mathbb{M}^f_g(t)\}_{t\ge0}$. Now we prove a version of the law of iterated logarithm for TCCPP. The following result will be used (see Theorem 14 of \cite{Bertoin2004}, p. 92).

\begin{lemma}\label{lil.Subordinator}
Let $\{D_f(t)\}_{t\ge 0}$ be a L\'{e}vy subordinator whose Laplace exponent $f$ is a regularly varying function at $0+$ with index $0<\gamma <1$, that is, $\lim_{x\rightarrow 0+}f(\lambda x)/f(x)=\lambda^{\gamma},~\lambda>0$. If
\begin{equation*}
h(t)=\frac{\log\log t}{\phi(t^{-1}\log\log t)}, \quad t>e
\end{equation*}
where $\phi$ is the inverse of $f$, then
\begin{equation*}
\lim\inf_{t\rightarrow\infty}\frac{D_f(t)}{h(t)} = \gamma(1-\gamma)^{(1-\gamma)/\gamma}\quad\text{a.s.} 
\end{equation*}
\end{lemma}
Now, we use Lemma \ref{lil.Subordinator} to obtain the law of iterated logarithm for $\{\mathbb{M}_f^g(t)\}_{t\ge 0}$.
\begin{theorem}\label{lil_TCCPP}
    The TCCPP $\{\mathbb{M}_f^g(t)\}_{t\ge 0}$ satisfies the following law of iterated logarithm:
    \begin{equation*}
        \liminf_{t\to\infty}\frac{\mathbb{M}_f^g(t)}{h(Y_g(t))}=\lambda\,\mathbb{E}[X_1]\,\gamma(1-\gamma)^{(1-\gamma)/\gamma}\quad \text{a.s.}\quad \text{such that}\quad Y_g(t)>e.
    \end{equation*}
\end{theorem}
\begin{proof}
Since the inverse subordinator $\{Y_g(t)\}_{t\ge0}$ is non-decreasing and
$Y_g(t)\to\infty$ almost surely as $t\to\infty$,  it follows from Lemma~\ref{lil.Subordinator} that,
\begin{equation}\label{lil_TCCPP_eq1}
\liminf_{t\to\infty}\frac{D_f(Y_g(t))}{h(Y_g(t))}=\gamma(1-\gamma)^{(1-\gamma)/\gamma}\quad \text{a.s.}
\end{equation}
Now using \eqref{lil_TCCPP_eq1} along with \eqref{strong_law_eq1} and \eqref{strong_law_eq2}, we obtain
\begin{align}
    \liminf_{t\to\infty}\frac{\mathbb{M}_f^g(t)}{h(Y_g(t))}&=\liminf_{t\to\infty}\frac{N_f^g(t)}{D_f(Y_g(t))}\frac{D_f(Y_g(t))}{h(Y_g(t))}\frac{1}{N_f^g(t)}\sum_{i=1}^{N_f^g(t)} X_i=\lambda\,\mathbb{E}[X_1]\,
\gamma(1-\gamma)^{(1-\gamma)/\gamma}
\quad \text{a.s.}\nonumber
\end{align}
This completes the proof.
\end{proof}
\begin{remark}
According to \citet{WilliamF2008} (see (2.6), p.~290), a p.g.f. $h(s)$ is infinitely divisible if and only if
\[
\log \frac{h(s)}{h(0)}=\sum_{k=1}^{\infty} a_k s^k, \quad a_k \ge 0,\quad \sum_{k=1}^{\infty} a_k < \infty.
\]
However, the p.g.f. \eqref{p.g.f.TCCPP} does not satisfy the above criterion. Hence, it is not infinitely divisible, which implies $\{\mathbb{M}_f^g(t)\}_{t\ge 0}$ is not infinitely divisible. Alternatively, note that $\{N_f^g(t)\}_{t\ge 0}$ is not a L\`{e}vy process and hence not infinitely divisible, which implies $\{\mathbb{M}_f^g(t)\}_{t\ge 0}$ is also not infinitely divisible.  
\end{remark}
Next, we obtain the asymptotic distribution of a scaled version of the TCCPP. The two lemmas presented below will be useful to obtain this distribution. 
\begin{lemma}\label{lem:rel-stab}
Let $\{D_f(t)\}_{t \ge 0}$ be a L\'{e}vy subordinator with Bernstein function $f$ satisfying $\mathbb{E}(D_f(1)) < \infty$, and let $\{Y_\alpha(t)\}_{t \ge 0}$ be the inverse $\alpha$-stable subordinator with $\alpha \in (0,1)$, independent of $\{D_f(t)\}_{t\ge0}$. Let $\{N(t)\}_{t \ge 0}$ be a homogeneous Poisson process with rate $\lambda > 0$, independent of $\{D_f(t)\}_{t \ge 0}$ and $\{Y_\alpha(t)\}_{t \ge 0}$. 
Then
\begin{equation*}
\frac{N(D_f(Y_\alpha(t)))}{C_f\,t^\alpha}\;\xrightarrow{P}\; 1\quad \text{as } t \to \infty, \quad \text{where}\quad C_f \;:=\; \frac{\lambda\,\mathbb{E}(D_f(1))}{\Gamma(1+\alpha)}.
\end{equation*}
\end{lemma}
\begin{proof}
Since $\{N(t)\}_{t \ge 0}$ is a Poisson process with rate $\lambda$ and $D_f(Y_\alpha(t))\xrightarrow{a.s.}\infty$ as $t\to\infty$, the weak law of large numbers gives
\begin{equation}
    \frac{N(D_f(Y_\alpha(t)))}{\lambda\,D_f(Y_\alpha(t))} \;\xrightarrow{p}\; 1
    \quad \text{as } t \to \infty.\label{lemma_term-1}
\end{equation}
Moreover since $\{D_f(t)\}_{t\ge0}$ is a L\'{e}vy process satisfying $\mathbb{E}(D_f(1))<\infty$ and $Y_\alpha(t)\overset{a.s.}{\longrightarrow}\infty$ as $t\to\infty$, the strong law of large numbers yields
\begin{equation*}
    \frac{D_f(Y_\alpha(t))}{Y_\alpha(t)} \;\xrightarrow{a.s.}\; \mathbb{E}(D_f(1)) \quad \text{as } t \to \infty,
\end{equation*}
which implies
\begin{equation}
    \frac{D_f(Y_\alpha(t))}{\mathbb{E}(D_f(1))\,Y_\alpha(t)} \;\xrightarrow{p}\; 1
    \quad \text{as } t \to \infty.\label{lemma_term-2}
\end{equation}
From Lemma 3 of \citet{Leonenko2019}, it follows that
\begin{equation}
    \frac{N(Y_\alpha(t))}{\mathbb{E}[N(Y_\alpha(t))]}
    \;=\;
    \frac{N(Y_\alpha(t))}{\lambda t^\alpha/\Gamma(1+\alpha)}
    \;\xrightarrow{p}\; 1
    \quad \text{as } t \to \infty.\label{lemma_term-2.5}
\end{equation}
Replacing $D_f(Y_{\alpha}(t))$ by $Y_{\alpha}(t)$ in \eqref{lemma_term-1} and using \eqref{lemma_term-2.5}, we have
\begin{equation}
    \frac{N(Y_\alpha(t))}{\lambda t^\alpha/\Gamma(1+\alpha)}\times\frac{\lambda Y_\alpha(t)}{N(Y_\alpha(t))}=\frac{\Gamma(1+\alpha)\,Y_\alpha(t)}{t^\alpha}=\frac{\mathbb{E}(D_f(1))\,Y_\alpha(t)}{(C_f/\lambda)\,t^\alpha} \;\xrightarrow{p}\; 1
    \quad \text{as } t \to \infty.\label{lemma_term-3}
\end{equation}
Now from \eqref{lemma_term-1}, \eqref{lemma_term-2} and \eqref{lemma_term-3}, we obtain
\begin{align*}
    \underbrace{\frac{N(D_f(Y_\alpha(t)))}{\lambda\,D_f(Y_\alpha(t))}}_{\xrightarrow{p}\,1}
    \times
    \underbrace{\frac{\lambda\,D_f(Y_\alpha(t))}{\lambda\,\mathbb{E}(D_f(1))\,Y_\alpha(t)}}_{\xrightarrow{p}\,1}
    \times
    \underbrace{\frac{\lambda\,\mathbb{E}(D_f(1))\,Y_\alpha(t)}{C_f\,t^\alpha}}_{\xrightarrow{p}\,1}=\frac{N(D_f(Y_\alpha(t)))}{C_f\,t^\alpha}\xrightarrow{p}\,1.
\end{align*}
\end{proof}
\begin{lemma}
    \label{lem:anscombe}
Let $\{X_i\}_{i \ge 1}$ be i.i.d. random variables in the domain of attraction of a stable law $\mu$, that is, there exist sequences
$(a_n)_{n \ge 1}$ and $(b_n)_{n \ge 1}$ of real constants such that
\begin{equation*}
    \bar{S}_n \;:=\; a_n S_n - b_n \;\xrightarrow{d}\; S \;\sim\; \mu,
\end{equation*} 
where $S_n = \sum_{k=1}^n X_k$. Then
\begin{equation*}
    a_{N(D_f(Y_\alpha(t)))}\,S_{N(D_f(Y_\alpha(t)))} - b_{N(D_f(Y_\alpha(t)))}
    \;\xrightarrow{d}\; S
    \quad \text{as } t \to \infty.
\end{equation*}
\end{lemma}
\begin{proof}
Consider the following three statements. \\ 
\noindent\textbf{(i)} The partial sums $S_n$ satisfy
\begin{equation*}
    \bar{S}_n = a_n S_n - b_n \;\xrightarrow{d}\; S \sim \mu
    \quad \text{as }\quad n \to \infty
\end{equation*}                           
by assumption, which satisfies first condition of Theorem~2 of \citet{Leonenko2019} (originally due to \citet{Anscombe1952}).

\noindent\textbf{(ii)} By Lemma~\ref{lem:rel-stab},
\begin{equation*}
    \frac{\widetilde{N}(t)}{C_f\,t^\alpha} \;\xrightarrow{P}\; 1
    \quad \text{as }\quad t \to \infty,
\end{equation*}
where $\widetilde{N}(t) := N(D_f(Y_\alpha(t)))$ and $\psi(t) := C_f\,t^\alpha \to \infty$ as $t \to \infty$. Hence $\widetilde{N}(t)$
is relatively stable with scaling function $\psi(t)$.

\noindent\textbf{(iii)} Since conditions \textbf{(i)} and \textbf{(ii)} hold for $(\bar{S}_n)_{n \ge 1}$, uniform continuity in probability of $(\bar{S}_n)_{n \ge 1}$ follows from Theorem~3 of \citet{Anscombe1952}, i.e.\ for every $\varepsilon > 0$ and $\eta > 0$ there exist $c = c(\varepsilon, \eta)$ and
$l_0 = l_0(\varepsilon, \eta)$ such that for all $l \ge l_0$,
\begin{equation*}
    \mathbb{P}(
        \max_{m\,:\,|m - l| < cl}
        |\bar{S}_m - \bar{S}_l| > \varepsilon
    ) < \eta.
\end{equation*}
Since all three conditions of Theorem~2 of \citet{Leonenko2019} are satisfied with $R_n := \bar{S}_n$ and random index $\widetilde{N}(t)$, we have
\begin{equation*}
    \bar{S}_{\widetilde{N}(t)}=a_{\widetilde{N}(t)}\,S_{\widetilde{N}(t)} - b_{\widetilde{N}(t)}
    \;\xrightarrow{d}\; S
    \quad \text{as }\quad t \to \infty.
\end{equation*}
\end{proof}
We now show that a suitably scaled process $\{\mathbb{M}_f^\alpha(t)=\mathbb{M}(D_f(Y_{\alpha}(t)))\}_{t\ge 0}$ (a particular case of the TCCPP) converges in distribution to a stable random variable.
\begin{theorem}\label{thm:stable-limit}
Let $\{X_i\}_{i\ge1}$ be a sequence of i.i.d.\ random variables with $\mathbb{E}(X_1)=0$. Let $S_n=\sum_{i=1}^nX_i$ and define $\widetilde{N}(t) := N(D_f(Y_\alpha(t)))$. Suppose there exist positive constants $(d_n)_{n\ge 1}$ such that
\begin{equation*}
    \frac{S_n}{d_n}\xrightarrow{d}S,
\end{equation*}
where $S$ is a stable random variable with stability index $\gamma\in(1,2]$. Then, as $t\to\infty$,
\begin{equation*}
    \frac{\mathbb{M}_f^\alpha(t)}{d_{\lfloor C_f t^\alpha \rfloor}}
    = \frac{1}{d_{\lfloor C_f t^\alpha \rfloor}}\sum_{i=1}^{\widetilde{N}(t)}X_i
    \xrightarrow{d}S,\quad \text{where}\quad C_f \;:=\; \frac{\lambda\,\mathbb{E}(D_f(1))}{\Gamma(1+\alpha)}.
\end{equation*}
\end{theorem}
\begin{proof} Since $b_n=0$ in $\tilde{S}_n$, Lemma~\ref{lem:anscombe} implies
\begin{equation}
\label{eq:random-norming}
    \frac{S_{\widetilde{N}(t)}}{d_{\widetilde{N}(t)}}=\frac{\sum_{i=1}^{\widetilde{N}(t)}X_i}{d_{\widetilde{N}(t)}}=\frac{\mathbb{M}_f^\alpha(t)}{d_{\widetilde{N}(t)}}
    \;\xrightarrow{d}\; S
    \quad \text{as } t \to \infty.
\end{equation}
Next consider Lemma~2.9 (a) of \citet{Gut1974}. Comparing Gut's notations with ours, we have $B_n := d_n$, $N := \widetilde{N}(t)$, and $B(c/\theta) := d_{\lfloor C_f t^\alpha \rfloor}$.
Since $B_N/B(c/\theta) \xrightarrow{p} 1$ by the above lemma, it follows that
\begin{equation}
\label{eq:ratio}
    \frac{d_{\widetilde{N}(t)}}{d_{\lfloor C_f\,t^\alpha \rfloor}}
    \;\xrightarrow{p}\; 1
    \quad \text{as } t \to \infty.
\end{equation}
Using Slutsky's theorem, from \eqref{eq:random-norming} and \eqref{eq:ratio}, we obtain
\begin{equation*}
    \frac{\mathbb{M}_f^\alpha(t)}{d_{\lfloor C_f\,t^\alpha \rfloor}}=\underbrace{\frac{d_{\widetilde{N}(t)}}{d_{\lfloor C_f\,t^\alpha \rfloor}}}_{\xrightarrow{p}\;1} \times \underbrace{\frac{\mathbb{M}_f^\alpha(t)}{d_{\widetilde{N}(t)}}
    }_{\xrightarrow{d}\;S}\xrightarrow{d} 1\times S \;=\; S
    \quad \text{as}\quad t \to \infty.
\end{equation*}
\end{proof}
\begin{remark}\leavevmode
\begin{itemize}
\item[(i)] The restriction $\gamma \in (1,2]$ ensures that the mean $\mathbb{E}(X_1)$ exists. It is assumed that the norming sequence $(d_n)$ is regularly varying, which is essential for replacing the random index $\widetilde{N}(t)$ by the deterministic sequence $\lfloor C_f t^\alpha \rfloor$. 

\item[(ii)] Theorem \ref{thm:stable-limit} does not require $\{X_i\}_{i\ge1}$ 
    to be independent of the process $\{N(D_f(Y_{\alpha}(t)))\}_{t\ge0}$. Moreover, the above result remains valid even for mixing sequences $X_1,X_2,\ldots$ instead of i.i.d. sequence (see \citet{Csorgo1973} for a generalization of Anscombe's theorem for mixing sequences).
\end{itemize}
\end{remark}
The results of Lemma~\ref{lem:rel-stab}, Lemma~\ref{lem:anscombe}, and 
Theorem \ref{thm:stable-limit} extend naturally to the TCCPP 
$\{\mathbb{M}_f^g(t)\}_{t\ge0}$.
Next, let $\mathcal{F}_0:=\sigma(D_f(Y_g(t)),t\ge 0)$, $\mathcal{F}^{N_f^g(t)}_t:=\sigma(N_f^g(s),\,0\leq s\leq t)$ and $\mathcal{F}_t=\mathcal{F}_0\vee \mathcal{F}^{N_f^g(t)}_t$.
The following result provides the asymptotic distribution of $\{N_f^g(t)\}_{t\ge 0}$.
\begin{proposition}
    If $\{N_f^g(t)\}_{t\ge 0}$ is a time changed Poisson process adapted to $\mathcal{F}_t$, then
    \begin{equation*}
\frac{N_f^g(t)-\lambda D_f(Y_g(t))}{\sqrt{\lambda D_f(Y_g(t))}}\stackrel{d}{\longrightarrow}  \mathcal{N}(0,1) \quad \text{as} \quad t \to\infty.
\end{equation*}
\end{proposition}
\begin{proof} Let $\{N(t)\}_{t\ge 0}$ be a Poisson process with rate $\lambda$ independent of $\{D_f(Y_g(t))\}_{t\ge 0}$. Define $X(t):=N_f^g(t)-\lambda D_f(Y_g(t))$.
First we show that $\{X(t)\}_{t\ge 0}$ is a martingale with respect to $\mathcal{F}_t$. For $0\le s<t$, using the tower property, we have
\begin{equation*}
    E(X(t)\mid\mathcal F_s)
=E(\,E(N_f^g(t)-\lambda D_f(Y_g(t))\;|\;\mathcal{F}_0\vee\mathcal F_s)\;|\;\mathcal F_s).
\end{equation*}
Note that $N_f^g(t)|D_f(Y_g(t))\sim \text{Poi}(\lambda D_f(Y_g(t)))$. Using the martingale property of a  Poisson process, we get
\begin{equation*}
    E[N_f^g(t)-\lambda D_f(Y_g(t))\;|\;\mathcal{F}_0\vee\mathcal F_s]
= N_f^g(s)-\lambda D_f(Y_g(s)).
\end{equation*}
Hence
\begin{equation*}
    E[X(t)\mid\mathcal F_s]=N_f^g(s)-\lambda D_f(Y_g(s))=X(s),
\end{equation*}
showing that $\{X(t)\}_{t\ge 0}$ is a martingale. It follows that $\{A(t)=\lambda D_f(Y_g(t))\}_{t\ge 0}$ is the compensator of $\{N_f^g(t)\}_{t\ge0}$, which is continuous in $t$. Now, let $f(u)\equiv1$.
Then
\begin{equation*}
B_t^2=\int_0^t [f(u)]^2\,dA(u)=\lambda D_f(Y_g(t))>0
\end{equation*}
and
\begin{equation*}
X_t:=\int_0^tf(u)\,[dN_f^g(u)-dA(u)]=N_f^g(t)-\lambda D_f(Y_g(t)).
\end{equation*}
Therefore
\begin{equation*}
\frac{X_t}{B_t}=\frac{N_f^g(t)-\lambda D_f(Y_g(t))}{\sqrt{\lambda D_f(Y_g(t))}}.
\end{equation*}
Applying Lemma 1 of \citet{Leonenko2019}, the randomly normed integral $X_t/B_t$ converges $\mathcal{F}_0$-stably (see Definition 2 of \citet{Leonenko2019}) as $t\, \to\infty$. Since $\mathcal{F}$-stable convergence implies weak convergence or convergence in distribution, we have
\begin{equation*}
\frac{N_f^g(t)-\lambda D_f(Y_g(t))}{\sqrt{\lambda D_f(Y_g(t))}}\stackrel{d}{\longrightarrow} \mathcal{N}(0,1)\quad \text{as} \quad t \to\infty.\label{convergence_TC_Poisson}
\end{equation*}    
\end{proof}
We now establish a central limit theorem for the TCCPP under suitable conditions.
\begin{theorem} 
Assume the following conditions hold:
\begin{equation*}
\frac{D_f(Y_g(t))-\mathbb{E}(D_f(Y_g(t)))}{\sqrt{\mathbb{V}(D_f(Y_g(t)))}}\xrightarrow{d}\mathcal N(0,1),\quad\frac{D_f(Y_g(t))}{\mathbb{E}(D_f(Y_g(t)))}\xrightarrow{p}1, \quad \frac{\mathbb{E}(D_f(Y_g(t)))}{[\mathbb{V}(N_f^g(t))]^{k/2}}\to \eta, 
\end{equation*}
where $\eta>0$ if $k=2$ or $\eta=0$ if $k\ge 3$ as $t\to\infty$. Then 
\begin{equation*}
    \frac{\mathbb{M}_f^g(t)-\mathbb{E}(\mathbb{M}_f^g(t))}{\sqrt{\mathbb{V}(\mathbb{M}_f^g(t))}}\xrightarrow{d}\mathcal N(0,1)\quad \text{as}\quad t\to\infty.
\end{equation*}
\end{theorem}
\begin{proof} 
Let $\mu_t=\mathbb{E}(D_f(Y_g(t)))<\infty$ and $\nu_t^2=\mathbb{V}(D_f(Y_g(t)))<\infty$. Using Lemma 2.1 of \citet{Leonenko2014}, we have $\mathbb{E}(N_f^g(t)=\lambda\mu_t<\infty$ and $\mathbb{V}(N_f^g(t))=\sigma^2_t=\lambda\mu_t+\lambda^2\nu_t^2<\infty$. Define 
\begin{equation*}
Y_t:=\frac{N_f^g(t)-\mathbb{E}(N_f^g(t))}{\sqrt{\mathbb{V}(N_f^g(t))}}=\frac{N_f^g(t)-\lambda\mu_t}{\sigma_t}.
\end{equation*}
The characteristic function of $Y_t$ is given by
\begin{align}
\phi_t(s)&=\mathbb{E}(\mathbb{E}(\exp(is\frac{N(D_f(Y_g(t))-\lambda\mu_t}{\sigma_t})|D_f(Y_g(t))))
=\mathbb{E}(\exp(\lambda D_f(Y_g(t))(e^{\frac{is}{\sigma_t}}-1)-is\frac{\lambda\mu_t}{\sigma_t}))\nonumber
\end{align}
where $i=\sqrt{-1}$. Let $\displaystyle R(t)=\sum_{k=3}^{\infty}\frac{(is)^k}{k!\sigma_t^k}$. Using Taylor series expansion, we get
\begin{align}
    \phi_t(s)&=\mathbb{E}(\exp(\lambda D_f(Y_g(t))(\frac{is}{\sigma_t}-\frac{s^2}{2\sigma^2_t}+R(t))-is\frac{\lambda\mu_t}{\sigma_t}))\nonumber\\
    &=\mathbb{E}(\exp(is\frac{\lambda (D_f(Y_g(t))-\mu_t)}{\sigma_t}-\frac{s^2}{2}\frac{\lambda D_f(Y_g(t))}{\sigma_t^2}+ \lambda D_f(Y_g(t)) R(t)))\nonumber\\
    &=\mathbb{E}(\exp(is\frac{(D_f(Y_g(t))-\mu_t)}{\nu_t}\sqrt{\frac{\lambda^2\nu_t^2}{\sigma_t^2}}-\frac{s^2}{2}\frac{D_f(Y_g(t))}{\mu_t}\frac{\lambda\mu_t}{\sigma_t^2}+ \lambda D_f(Y_g(t)) R(t)))\label{char_Yt_eq1}.
\end{align}
By assumption, $\mu_t/\sigma^k_t\to 0$ for $k\ge 3$ which implies
\begin{equation}
\mathbb{E}(\lambda D_f(Y_g(t)) R(t)) = \lambda\mu_tR(t) = \lambda\mu_t\sum_{k=3}^{\infty}\frac{(is)^k}{k!\sigma_t^k}\to 0\implies\lambda D_f(Y_g(t))R(t)\xrightarrow{p}0.\label{char_fun_third_term}
\end{equation}
Next, we have
\begin{equation}
\frac{D_f(Y_g(t))}{\mu_t}\xrightarrow{p}1\implies -\frac{s^2}{2}\frac{D_f(Y_g(t))}{\mu_t}\frac{\lambda\mu_t}{\sigma_t^2}\xrightarrow{p}-\frac{s^2}{2}\frac{\lambda\mu_t}{\sigma_t^2}.\label{char_fun_second_term}
\end{equation}
From \eqref{char_fun_third_term} and \eqref{char_fun_second_term}, it follows that
\begin{equation*}
-\frac{s^2}{2}\frac{D_f(Y_g(t))}{\mu_t}\frac{\lambda\mu_t}{\sigma_t^2}+\lambda D_f(Y_g(t))R(t)\xrightarrow{p}-\frac{s^2}{2}\frac{\lambda\mu_t}{\sigma_t^2}.
\end{equation*}
By assumption, $\mu_t/\sigma^2_t\to\eta$ which implies
\begin{equation}\label{char_fun_2+3}
    -\frac{s^2}{2}\frac{D_f(Y_g(t))}{\mu_t}\frac{\lambda\mu_t}{\sigma_t^2}+\lambda D_f(Y_g(t))R(t)\xrightarrow{p}\frac{-s^2\lambda\eta}{2}.
\end{equation} 
Now let $Z\sim \mathcal{N}(0,1)$. Then 
\begin{equation*}
\frac{D_f(Y_g(t))-\mu_t}{\nu_t}\xrightarrow{d}Z\implies is\frac{D_f(Y_g(t))-\mu_t}{\nu_t}\sqrt{\frac{\lambda^2\nu_t^2}
{\sigma_t^2}}\xrightarrow{d}isZ\sqrt{\frac{\lambda^2\nu_t^2}{\sigma_t^2}}.
\end{equation*}
Since $\mu_t/\sigma^2_t\to\eta$ by assumption, we have $1-\lambda\mu_t/\sigma^2_t=\lambda^2\nu^2_t/\sigma^2_t\to 1-\lambda\eta$, which implies
\begin{equation}
is\frac{D_f(Y_g(t))-\mu_t}{\nu_t}\sqrt{\frac{\lambda^2\nu_t^2}
{\sigma_t^2}}\xrightarrow{d}isZ\sqrt{1-\lambda\eta}.\label{char_fun_1}
\end{equation}
By Slutsky's theorem, it follows from \eqref{char_fun_2+3} and \eqref{char_fun_1} that
\begin{equation}
is\frac{D_f(Y_g(t))-\mu_t}{\nu_t}\sqrt{\frac{\lambda^2\nu_t^2}{\sigma_t^2}}-\frac{s^2}{2}\frac{D_f(Y_g(t))}{\mu_t}\frac{\lambda\mu_t}{\sigma_t^2}+\lambda D_f(Y_g(t))R(t)\xrightarrow{d}isZ\sqrt{1-\lambda\eta}-\frac{s^2\lambda\eta}{2}.\nonumber
\end{equation}
Using Portmanteau theorem, we have
\begin{equation}
\lim_{t\to\infty}\phi_t(s)
=\mathbb{E}(\exp(isZ\sqrt{1-\lambda\eta}-\frac{s^2\lambda\eta}{2}))=\exp(\frac{-s^2(1-\lambda\eta)}{2})\exp(\frac{-s^2\lambda\eta}{2})=\exp(\frac{-s^2}{2})\label{char_fun_Yt}
\end{equation}
showing that $N_f^g(t)\sim\mathcal{N}(\mu_t,\,\sigma_t^2)$ asymptotically. Next, note that the distribution function of $N_f^g(t)$ is 
\begin{equation*}
\mathbb{P}(N_f^g(t)\le x) = \int_{0}^{\infty}\mathbb{P}(N(u)\le x)f_{D_f(Y_g(t))}(u)du = \int_{0}^{\infty}\sum_{k=0}^{x}\frac{e^{-\lambda u}(\lambda u)^k}{k!}f_{D_f(Y_g(t))}(u)du
\end{equation*}
which depends on $\lambda$. Moreover
\begin{align}
    \frac{\sigma_t}{\mathbb{V}(\mathbb{M}_f^g(t))}=\frac{\sqrt{\lambda\mu_t+\lambda^2\nu_t^2}}{\lambda\mathbb{E}(X_1^2)\mu_t+\lambda^2(\mathbb{E}(X_1))^2\nu_t^2}=&\frac{1}{\lambda\nu_t}\frac{\sqrt{1+\frac{\mu_t}{\lambda\nu_t^2}}}{((\mathbb{E}(X_1))^2+\frac{\mathbb{E}(X_1^2)\mu_t}{\lambda\nu_t^2})}\rightarrow 0 \quad\text{as}\quad \lambda\to\infty,\label{order_sigma_var}
\end{align}
that is, $\sigma_t=o(\mathbb{V}(\mathbb{M}_f^g(t)))$ and
\begin{equation}    \mathbb{V}(\mathbb{M}_f^g(t))=\lambda\mathbb{E}(X_1^2)\mu_t+\lambda^2(\mathbb{E}(X_1))^2\nu_t^2\to \infty\quad\text{as}\quad\lambda\to\infty.\label{var_to_infty}
\end{equation}
From \eqref{char_fun_Yt}, \eqref{order_sigma_var} and \eqref{var_to_infty}, we observe that the assumptions in Corollary $2$ of \citet{Robbins1948} are satisfied and hence $\{\mathbb{M}_f^g(t)\}_{t\ge 0}$ is asymptotically normally distributed with mean $\mathbb{E}(\mathbb{M}_f^g(t))$ and variance $\mathbb{V}(\mathbb{M}_f^g(t))$. Alternatively, the asymptotic normality of $\{N(D_f(Y_g(t))\}_{t\ge 0}$ (see \eqref{char_fun_Yt}) implies the asymptotic normality of $\{\mathbb{M}_f^g(t)\}_{t\ge 0}$ using Corollary 4 of \citet{Robbins1948}. Therefore 
\begin{equation*}
    \frac{\mathbb{M}_f^g(t)-\mathbb{E}(\mathbb{M}_f^g(t))}{\sqrt{\mathbb{V}(\mathbb{M}_f^g(t))}}\xrightarrow{d}\mathcal N(0,1)\quad \text{as}\quad t\to\infty.
\end{equation*} 
\end{proof}

\section{Characterizations}\label{sec-6}
In this section, we establish a martingale characterization of the TCCPP, followed by characterizations based on a univariate multinomial process, a triangular array process, and the waiting time distribution.
\subsection{Martingale characterization}
We first discuss the martingale characterization of the classical CPP as established in \citet{Zhang2016} and provide an alternative proof.
\begin{theorem}\label{martingale_char_CPP}
The process $\{\mathbb{M}(t)\}_{t\ge 0}$ is a CPP if and only if $\{\mathbb{M}(t)-\lambda t \mathbb{E}(X_1)\}_{t\ge 0}$ is a $\mathcal{F}_t$-martingale where $\mathcal{F}_t$ is the natural filtration of $\{\mathbb{M}(t)\}_{t\ge 0}$.
\end{theorem}
\begin{proof}
Let $\{\mathbb{M}(t)\}_{t\ge 0}$ be a CPP. For $s\le t$, consider
\begin{align}
\mathbb{E}(\mathbb{M}(t)-\lambda t\mathbb{E}(X_1)|\mathcal{F}_s)&=\mathbb{E}(\mathbb{M}(s)-\lambda s\mathbb{E}(X_1)|\mathcal{F}_s)\nonumber+\mathbb{E}(\mathbb{M}(t)-\mathbb{M}(s)|\mathcal{F}_s)-(\lambda (t-s)\mathbb{E}(X_1)|\mathcal{F}_s)\nonumber\\
&=\mathbb{M}(s)-\lambda s\mathbb{E}(X_1)+\lambda (t-s)\mathbb{E}(X_1)-\lambda (t-s)\mathbb{E}(X_1)\nonumber\\
&=\mathbb{M}(s)-\lambda s\mathbb{E}(X_1).\nonumber
\end{align}
Thus $\{\mathbb{M}(t)-\lambda t\mathbb{E}(X_1)\}_{t\ge 0}$ is a $\mathcal{F}_t$- martingale. 

Conversely suppose that $\{\mathbb{M}(t)-\lambda t\mathbb{E}(X_1)\}_{t\ge 0}$ is a $\mathcal{F}_t$-martingale. Define $\{Z(t)=\exp\{u\mathbb{M}(t)-\lambda t[\mathbb{E}(e^{uX_1})-1]\}\}_{t\ge 0}$. We have
\begin{align}
\mathbb{E}(Z(t)|\mathcal{F}_s)&=\mathbb{E}(Z(s)\exp\{u[\mathbb{M}(t)-\mathbb{M}(s)]-\lambda(t-s)(\mathbb{E}(e^{uX_1})-1)\}|\mathcal{F}_s)\nonumber\\
&=\mathbb{E}( \exp\{u\mathbb{M}(s)-\lambda s[\mathbb{E}(e^{uX_1})-1]\}\times\exp\{u[\mathbb{M}(t)-\mathbb{M}(s)]-\lambda (t-s)[\mathbb{E}(e^{uX_1})-1]\}|\mathcal{F}_s)\nonumber\\
&=Z(s)\mathbb{E}(\exp\{u[\mathbb{M}(t)-\mathbb{M}(s)]\}|\mathcal{F}_s)\mathbb{E}(\exp\{-\lambda(t-s)(\mathbb{E}(e^{uX_1})-1)\}|\mathcal{F}_s)\nonumber\\
&=Z(s)\mathbb{E}(\mathbb{E}(\exp\{u[\mathbb{M}(t)-\mathbb{M}(s)]\}|\mathcal{F}_s,N(t)-N(s)=k)|\mathcal{F}_s)\exp\{-\lambda(t-s)(\mathbb{E}(e^{uX_1})-1)\}\nonumber\\
&\hspace{9cm}(\because \mathcal{F}_s\subseteq\sigma(\mathcal{F}_s,N(t)-N(s)))\nonumber\\
&=Z(s)\mathbb{E}((\mathbb{E}(e^{uX_1}))^k|\mathcal{F}_s)\exp\{-\lambda(t-s)(\mathbb{E}(e^{uX_1})-1)\}\nonumber\\
&=Z(s)\exp\{\lambda(t-s)(\mathbb{E}(e^{uX_1})-1)\}\exp\{-\lambda(t-s)(\mathbb{E}(e^{uX_1})-1)\}\nonumber\\
&\hspace{9cm}(\because N(t)-N(s))|\mathcal{F}_s\sim \text{Poi}(\lambda(t-s))\nonumber\\
&=Z(s).\nonumber
\end{align}
Thus $\{Z(t)=\exp\{u\mathbb{M}(t)-\lambda t[\mathbb{E}(e^{uX_1})-1]\}\}_{t\ge 0}$ is a martingale. By definition of a martingale, we have 
\begin{align}
\mathbb{E}(Z(t)|\mathcal{F}_s)&=Z(s)\nonumber\\
\implies\mathbb{E}(Z(s)\exp\{u[\mathbb{M}(t)-\mathbb{M}(s)]-\lambda(t-s)(\mathbb{E}(e^{uX_1})-1)\}|\mathcal{F}_s)&=Z(s)\nonumber\\
\implies\mathbb{E}(e^{u[\mathbb{M}(t)-\mathbb{M}(s)]}|\mathcal{F}_s)&=e^{\lambda(t-s)(\mathbb{E}(e^{uX_1})-1)}.\label{CPP_martingale_eq1}
\end{align}
Hence, $\mathbb{M}(t)-\mathbb{M}(s)\mid \mathcal{F}_s \sim \mathrm{Compound\ Poisson}(\lambda(t-s),f_X)$, since a compound Poisson random variable $Y$ with parameter $\gamma$ and jump distribution $f_x$ has moment generating function $\mathbb{E}[e^{uY}] = \exp\!(\gamma(\mathbb{E}[e^{uX_1}] - 1))$. Taking expectation on both sides of \eqref{CPP_martingale_eq1}, we obtain $\mathbb{M}(t)-\mathbb{M}(s) \sim \mathrm{Compound\ Poisson}(\lambda(t-s),f_X)$. Since the distribution of the increment depends only on $(t-s)$, the increments of $\{\mathbb{M}(t)\}_{t\ge0}$ are stationary. Moreover, as the conditional distribution of $\mathbb{M}(t)-\mathbb{M}(s)$ given $\mathcal{F}_s$ does not depend on $\mathcal{F}_s$, the increments are independent. Therefore, $\{\mathbb{M}(t)\}_{t\ge0}$ has stationary and independent increments with compound Poisson distributed increments, and hence $\{\mathbb{M}(t)\}_{t\ge0}$ is a CPP. This completes the proof.
\end{proof}
Next, we establish the martingale characterization of a TCCPP. We begin with the following lemma (see Lemma 1 of \citet{Aletti2018}) which will be used later.
\begin{lemma}\label{martingale_lemma}
    Let $X$ be a right-continuous martingale. If $T$ and $S$ are stopping times such that $\mathbb{P}(T<\infty)=1$ and $\{X(t\wedge T),t\ge 0\}$ is uniformly integrable, then $\mathbb{E}(X(T)|\mathcal{F}_{S\wedge T})=X(S\wedge T)$.
\end{lemma}
\begin{theorem}\label{martingale_char_TCCPP}
    Let $\{K(t)\}_{t\ge 0}$ be a pure jump process with finitely many jumps on compact intervals. Then $\{K(t)\}_{t\ge 0}$ is a TCCPP if and only if there exist a CPP $\{\mathbb{M}(t)\}_{t\ge0}$ with rate $\lambda>0$ and an independent subordinator $\{D_f(Y_g(t))\}_{t\ge0}$ such that the process \begin{equation}\label{martingale_form_tcpp} \{X(t)\}_{t\ge0} = \{K(t)-\lambda \mathbb{E}(X_1)D_f(Y_g(t))\}_{t\ge0} \end{equation} is a right–continuous martingale with respect to the induced filtration \[ \mathcal{F}_t=\sigma(K(s),s\le t)\vee\sigma(D_f(Y_g(s)),s\ge0) \] and for any $T>0$, \begin{equation}\label{stopping_times_martingale_tccpp} \{X(\tau),\,\tau\,\, \text{is stopping time such that}\,\,D_f(Y_g(\tau))\le T\} \end{equation} is uniformly integrable.
\end{theorem}
\begin{proof}
    Let $\{K(t)\}_{t\ge0}$ be a TCCPP such that $K(t)=\mathbb{M}(D_f(Y_g(t)))$, where $\{\mathbb{M}(t)\}_{t\ge0}$ is a CPP independent of the time change $\{D_f(Y_g(t))\}_{t\ge0}$. Then 
    \begin{equation*}
         X(t):=K(t)-\lambda \mathbb{E}(X_1)D_f(Y_g(t)) = \mathbb{M}(D_f(Y_g(t)))-\lambda \mathbb{E}(X_1)D_f(Y_g(t)).
    \end{equation*}
    First we show that $\{X(t)\}_{t\ge0}$ is a martingale with respect to the filtration $\mathcal{F}_t=\sigma(K(s),s\le t)\vee\sigma(D_f(Y_g(s)),s\ge0)$. For $0\le s<t$, using the tower property and conditioning on $D_f(Y_g(t))$, we have 
    \begin{equation}
        \mathbb{E}(X(t)\mid\mathcal{F}_s) = \mathbb{E}(\mathbb{E}(\mathbb{M}(D_f(Y_g(t)))-\lambda \mathbb{E}(X_1)D_f(Y_g(t))\,|\,D_f(Y_g(t))\vee\mathcal{F}_s)|\mathcal{F}_s). \nonumber
    \end{equation} 
    Note that conditional on $D_f(Y_g(t))$, the random variable $\mathbb{M}(D_f(Y_g(t)))$ has a compound Poisson distribution with parameter $\lambda D_f(Y_g(t))$. Since the compensated CPP $\{\mathbb{M}(t)-\lambda t\mathbb{E}(X_1)\}_{t\ge0}$ is a martingale, we obtain 
    \begin{equation}
         \mathbb{E}(\mathbb{M}(D_f(Y_g(t)))-\lambda \mathbb{E}(X_1)D_f(Y_g(t))\,|\,D_f(Y_g(t))\vee\mathcal{F}_s) = \mathbb{M}(D_f(Y_g(s)))-\lambda \mathbb{E}(X_1)D_f(Y_g(s)). \nonumber
    \end{equation}
    Therefore $\mathbb{E}[X(t)\mid\mathcal{F}_s] = \mathbb{M}(D_f(Y_g(s)))-\lambda \mathbb{E}(X_1)D_f(Y_g(s)) = X(s), $ which proves that $\{X(t)\}_{t\ge0}$ is a martingale. Right continuity follows from the right continuity of the CPP and the time change process. 

    Next we show that the family $\{X(\tau):D_f(Y_g(\tau))\le T\}$ is uniformly integrable. Let $\tau$ be a stopping time such that $D_f(Y_g(\tau))\le T$. Then 
    \begin{equation}\label{martingale_char_TCCPP_eq1}
        \mathbb{E}[X(\tau)] = \mathbb{E}[\mathbb{M}(D_f(Y_g(\tau)))-\lambda \mathbb{E}(X_1)D_f(Y_g(\tau))] = 0.
    \end{equation} 
    Using the law of total variance, 
    \begin{align}
        \mathbb{V}(X(\tau)) &= \mathbb{E}(\mathbb{V}(K(\tau)-\lambda \mathbb{E}(X_1)D_f(Y_g(\tau))\mid D_f(Y_g(\tau))))+\mathbb{V}(\mathbb{E}(K(\tau)-\lambda \mathbb{E}(X_1)D_f(Y_g(\tau))\mid D_f(Y_g(\tau))))\nonumber\\
        &=\lambda \mathbb{E}(X_1^2)\mathbb{E}(D_f(Y_g(\tau)))+0\nonumber
    \end{align}
    Therefore using \eqref{martingale_char_TCCPP_eq1}, we get $\displaystyle\mathbb{E}(X^2(\tau))=\lambda \mathbb{E}(X_1^2)\mathbb{E}(D_f(Y_g(\tau)))\le \lambda \mathbb{E}(X_1^2)T.$
    Hence the family \eqref{stopping_times_martingale_tccpp} is uniformly bounded in $L^2$ and therefore uniformly integrable.

    Conversely, suppose that the process in \eqref{martingale_form_tcpp} satisfies the martingale property and the family \eqref{stopping_times_martingale_tccpp} is uniformly integrable. It is enough to prove that $K(t)=\mathbb{M}(D_f(Y_g(t)))$, where $\{\mathbb{M}(t)\}_{t\ge0}$ is a CPP independent of the process $\{D_f(Y_g(t))\}_{t\ge0}$ with rate $\lambda>0$ and jump distribution $f_x$. Define $Z(t)=\inf\{s: D_f(Y_g(s))\ge t\}$, the inverse of $\{D_f(Y_g(t))\}_{t\ge0}$. Then $\{Z(t)\}_{t\ge0}$ forms a family of stopping times. Hence, by Lemma \ref{martingale_lemma}, the process $X(Z(t))=K(Z(t))-\lambda \mathbb{E}(X_1)D_f(Y_g(Z(t)))$ is also a martingale. Since the time-change process is non-decreasing and continuous, we have $D_f(Y_g(Z(t)))=t$ a.s., and therefore $\{K(Z(t))-\lambda \mathbb{E}(X_1)t\}_{t\ge0}$ is a martingale. Moreover, as $\{Z(t)\}_{t\ge0}$ is increasing, the process $\{K(Z(t))\}_{t\ge0}$ is a pure jump process with finitely many jumps on compact intervals. By Theorem \ref{martingale_char_CPP}, it follows that $\{K(Z(t))\}_{t\ge0}$ is a classical CPP with rate $\lambda$ and jump distribution $f_x$. If we define $\mathbb{M}(t)=K(Z(t))$, then $K(t)=\mathbb{M}(D_f(Y_g(t)))$, and hence $\{K(t)\}_{t\ge0}$ is a TCCPP. This completes the proof.
\end{proof}

\subsection{Characterization by univariate multinomial process}
Consider a sequence of \(n\) independent trials, each resulting in one of \((k+1)\) mutually exclusive and exhaustive outcomes with corresponding probabilities \(p_0(t), p_1(t), \dots, p_k(t)\) at time $t>0$ satisfying $\sum_{i=0}^{k} p_i(t) = 1$. Suppose the occurrence of $i^{th}$ outcome in a trial corresponds to $i$ successes for $i = 0,1,\dots,k$. Then the total number of successes \(X(t)\) in \(n\) trials follows a univariate multinomial distribution (see Section 11.2.18 of \citet{Johnson2005}) with parameters $n,k,p_0(t),p_1(t),\dots,p_k(t)$. The p.m.f. of $X(t)$ is given by
\begin{equation}
    \mathbb{P}(X(t)=x)=\sum_{\substack{\sum_{i=0}^{k}r_i=n\\\sum_{i=0}^{k}ir_i=x}}\frac{n!}{r_0!r_1!\cdots r_k!}(p_0(t))^{r_0}(p_1(t))^{r_1}\cdots(p_k(t))^{r_k},\quad x=0,1,\dots,nk.\notag
\end{equation}
where $r_0,r_1,\dots,r_k\in\mathbb{N}_0$ and $n,k\in\mathbb{N}$. We refer to the collection $\{ X(t) \}_{t \ge 0}$ as a univariate multinomial process. The following result gives the asymptotic distribution of a time-changed univariate multinomial process $\{ X(D_f(Y_g(t))) \}_{t \ge 0}$, where $\{X(t)\}_{t\ge 0}$ is assumed to be independent of $\{D_f(Y_g(t))\}_{t\ge 0}$. 
We show that the TCCPP arises as the limit of a sequence of univariate time-changed multinomial processes when the number of trials grows and outcome probabilities shrink.
\begin{theorem}
The time-changed univariate multinomial process $\{X(D_f(Y_g(t)))\}_{t \ge 0}$ converges in distribution to a TCCPP as $n \to \infty$ and $p_i(t) \to 0$ such that $\,n\,p_i(t) \to \lambda\,\alpha_i\,t$, where $\lambda$ is the intensity and $\mathbb{P}(X_1 = i) = \alpha_i$ is the jump size distribution of the TCCPP for $1\le i\le k$.
\end{theorem}
\begin{proof}
For $|z|\le 1$, the p.g.f. of the univariate multinomial process $\{ X(t) \}_{t \ge 0}$ is
\begin{align}
G(z,t)
= \sum_{x=0}^{nk} z^x \, \mathbb{P}(X(t)=x)= \sum_{\substack{\sum_{i=0}^{k} r_i = n \\ \sum_{i=0}^{k} i r_i = x}}
\frac{n!}{r_0! r_1! \cdots r_k!}
\prod_{i=0}^{k} (p_i(t))^{r_i}
z^{\sum_{i=0}^{k} i r_i}
&= \sum_{\substack{\sum_{i=0}^{k} r_i = n}}
\frac{n!}{r_0! r_1! \cdots r_k!}
\prod_{i=0}^{k} (p_i(t) z^i)^{r_i} \notag\\
&= ( \sum_{i=0}^{k} p_i(t) z^i )^n.\notag 
\end{align}
Since $\{X(t)\}_{t\ge 0}$ is independent of $\{D_f(Y_g(t)\}_{t\ge 0}$, the p.g.f. of the time-changed univariate multinomial process $\{X(D_f(Y_g(t)))\}_{t \ge 0}$ is given by
   \begin{align}
\Tilde{G}(z,t)=\int_{0}^{\infty}G(z,u)f_{D_f(Y_g(t))}(u)\,du= \int_{0}^{\infty}(\sum_{i=0}^{k} p_i(u)\, z^i)^n f_{D_f(Y_g(t))}(u)\,du.\label{multinomial_pgf_eq_2}
   \end{align}
We rewrite the first term in the above integral as
\begin{align}
\sum_{i=0}^{k} p_i(u) z^i
=
p_0(u) + \sum_{i=1}^{k} p_i(u) z^i
&=
1 - \sum_{i=1}^{k} p_i(u) + \sum_{i=1}^{k} p_i(u) z^i=
1 + \sum_{i=1}^{k} p_i(u) (z^i - 1).\label{multinomial_pgf_eq_3}
\end{align}
Substituting \eqref{multinomial_pgf_eq_3} into \eqref{multinomial_pgf_eq_2}, we have
\begin{align}
   \Tilde{G}(z,t)= \int_{0}^{\infty}(1 + \sum_{i=1}^{k} p_i(u) (z^i - 1))^n f_{D_f(Y_g(t))}(u)\,du. \label{multinomial_pgf_eq_4}
\end{align}

If $n \to \infty$ and $p_i(t) \to 0$ such that $n\,p_i(t) \to \lambda\,\alpha_i\,t$ for $1\le i\le k$, then
\[
\sum_{i=1}^{k} n\,p_i(u) (z^i - 1)
\;\longrightarrow\;
\lambda\,u\sum_{i=1}^{k} \alpha_i (z^i - 1).
\]
Therefore rewriting the integrand in \eqref{multinomial_pgf_eq_4} and using the standard result $\displaystyle\lim_{n \to \infty}
(1+\frac{x}{n})^n=e^x$, we have
\begin{align}
\lim_{n\to\infty}\Tilde{G}(z,t)
=\int_{0}^{\infty}(1 + \frac{1}{n}\lambda\,u\sum_{i=1}^{k} \alpha_i (z^i - 1))^n f_{D_f(Y_g(t))}(u)\,du&=\int_{0}^{\infty}\exp(\lambda\,u\sum_{i=1}^{k} \alpha_i (z^i - 1))\, f_{D_f(Y_g(t))}(u)\,du\notag\\
&=\mathbb{E}(\exp(\lambda\,D_f(Y_g(t))\sum_{i=1}^{k} \alpha_i (z^i - 1))).\label{p.g.f.limiting_process}
\end{align}
Now consider a TCCPP $\{\mathbb{M}_f^g(t)\}_{t\ge 0}$, where $\mathbb{M}_f^g(t)=\displaystyle\sum_{j=1}^{N_f^g(t)} X_j$, with intensity parameter $\lambda$ and jump distribution $\mathbb{P}(X_1=i)=\alpha_i$ for $1\le i\le k$. Using \eqref{m.g.f.TCCPP}, the p.g.f. of $\{\mathbb{M}_f^g(t)\}_{t\ge 0}$ is    
\begin{equation}
    \mathbb{E}(z^{\mathbb{M}_f^g(t)})=\mathbb{E}(\exp(\lambda D_f(Y_g(t))(\mathbb{E}(z^{X_1})-1)))=\mathbb{E}(\exp(\lambda\,D_f(Y_g(t))\sum_{i=1}^{k} \alpha_i (z^i - 1))).\label{p.g.f.TCCPP}
\end{equation}
Comparing \eqref{p.g.f.limiting_process} with \eqref{p.g.f.TCCPP}, we observe that the time-changed univariate multinomial process converges in distribution to the TCCPP $\{\mathbb{M}_f^g(t)\}_{t\ge 0}$. 
\end{proof}
\subsection{Characterization by triangular array process}
An infinite double array $[X_{ij}]$ of real-valued random variables is called a triangular array (see \citet{Cinlar2011}) if, for each $i\geq1$, there exists an integer $k_i$ such that
\begin{equation*}
    X_{ij}=0,\quad j>k_i,
\end{equation*}
where $\{k_i\}_{i\geq1}$ is an increasing sequence satisfying
$k_i\to\infty$ as $i\to\infty$. The random variables $X_{i1},X_{i2},\ldots,X_{ik_i}$ in the $i$th row are mutually independent. However, different rows may depend on each other.

As an extension of the above concept, we define the collection $\{X_{ij}(t)\}_{t\geq0},~i,j\geq1$ to be a triangular array process, where for each fixed
$t\geq0$, the random variables $X_{ij}(t)$ are elements of a triangular array. Define the $i^{th}$ row-sum process of $\{X_{ij}(t)\}_{t\geq0}$ by
\begin{equation*}
    Y_i(t):=\sum_{j=1}^{\infty}X_{ij}(t),\quad t\geq0.
\end{equation*}
In this section, we provide a characterization of the time-changed compound Poisson process in terms of a triangular array process. In particular, we show that the time-changed row-sum process $\{Y_i(D_f(Y_g(t)))\}_{t\ge 0}$ converges to the TCCPP for a fixed time $t\geq0$. 
\begin{theorem}\label{triangular_TCCPP}
    Let $\{X_{ij}(t)\}_{t\geq0}$ be a triangular array process of non-negative discrete random variables and let $\{Y_i(t)\}_{t\ge 0}$ be its $i^{th}$ row-sum process. Suppose for fixed $t\geq0$,
\begin{align}
    &\lim_{i\to\infty}\max_j \mathbb{P}(X_{ij}(t)\geq1)=0, \label{eq:TCCPP-triangular-1}\\
    &\lim_{i\to\infty}\sum_{j=1}^{\infty}\mathbb{P}(X_{ij}(t)\geq r+1\ge 2)=0,\quad r\in \mathbb{N}, \label{eq:TCCPP-triangular-2}\\
    &\lim_{i\to\infty}\sum_{j=1}^{\infty}\mathbb{P}(X_{ij}(t)=k)=\lambda t p_k,\quad k=1,\ldots,r. \label{eq:TCCPP-triangular-3}
\end{align}
If $D_f(Y_g(t))$ is independent of $Y_i(t)$, then $Y_i(D_f(Y_g(t)))
    \xrightarrow{d}\mathbb{M}_f^g(t)$, a TCCPP with intensity $\lambda$ and jump size distribution $\mathbb{P}(X_j=k)=p_k$, as $i\to\infty$.
\end{theorem}
\begin{proof}
    Let $B_{ij}(t)=\min\{X_{ij}(t),r\}$ and $S_i(t)=\sum_{j=1}^{\infty}B_{ij}(t)$. Then we have $\{Y_i(t)\neq S_i(t)\}
    \subseteq \bigcup_{j=1}^{\infty}
    \{X_{ij}(t)\neq B_{ij}(t)\}=
    \bigcup_{j=1}^{\infty}
    \{X_{ij}(t)\geq r+1\}.$ Consequently, as $i\to\infty$,
\begin{align}
    P\{Y_i(t)\neq S_i(t)\}
    &\leq
    \sum_{j=1}^{\infty}
    P\{X_{ij}(t)\geq r+1\}
    \to 0. \quad (\text{using}\,\eqref{eq:TCCPP-triangular-2})
    \label{eq:TCCPP-truncation}
\end{align}
Hence, $Y_i(t)-S_i(t)\xrightarrow{p}0$ so $Y_i(t)$ and $S_i(t)$ have the same limiting distribution. Therefore, it suffices to establish the limiting distribution of $S_i(t)$. Since the random variables in each row are mutually independent, we have
\begin{align}
    \mathbb{E}[z^{S_i(t)}]=\prod_{j=1}^{\infty}
    \mathbb{E}[z^{B_{ij}(t)}]&=\prod_{j=1}^{\infty}[\sum_{k=0}^{r}\mathbb{P}\{B_{ij}(t)=k\}z^k].
    \label{eq:TCCPP-pgf-product}
\end{align}
For $|z|\leq1$, define $x_{ij}(t)=\sum_{k=1}^{r}\mathbb{P}\{X_{ij}(t)=k\}(1-z^k)$. Since $\mathbb{P}\{B_{ij}(t)=0\}=1-\sum_{k=1}^{r}P\{X_{ij}(t)=k\}$, the expression \eqref{eq:TCCPP-pgf-product} can be written as
\begin{equation*}
\mathbb{E}[z^{S_i(t)}]=\prod_{j=1}^{\infty}(1-x_{ij}(t)).
\end{equation*}
By \eqref{eq:TCCPP-triangular-1},
\begin{align}
    \lim_{i\to\infty}\max_j x_{ij}(t)=\sum_{k=1}^{r}\lim_{i\to\infty}\max_j\mathbb{P}\{X_{ij}(t)=k\}(1-z^k)&\leq \sum_{k=1}^{r}
    (1-z^k)\lim_{i\to\infty}\max_j P\{X_{ij}(t)\geq1\}=0.
    \label{eq:TCCPP-x-infinitesimal}
\end{align}
Thus, for sufficiently large $i$, $x_{ij}(t)\leq 1/2$ for all $j$. 
By Taylor's series expansion, $\log(1-x)+x = -\sum_{k=2}^\infty x^k/k$. Hence 
$|\log(1-x)+x| \leq x^2\sum_{k=0}^\infty x^k = x^2/(1-x) \leq 2x^2$ for $0\leq x < 1/2$. Hence using \eqref{eq:TCCPP-x-infinitesimal}, we have 
\begin{equation*}
    |\sum_{j=1}^{\infty}x_{ij}(t)+\sum_{j=1}^{\infty}\log(1-x_{ij}(t))|\leq2\sum_{j=1}^{\infty}x_{ij}(t)^2\le 2(\max_jx_{ij}(t))
    \sum_{j=1}^{\infty}x_{ij}(t)
    \longrightarrow0.
\end{equation*}
Hence, 
\begin{align}
    \sum_{j=1}^{\infty}\log(1-x_{ij}(t))=-\sum_{j=1}^{\infty}x_{ij}(t)+o(1)&=-\sum_{k=1}^{r}(1-z^k)
    \sum_{j=1}^{\infty}
    \mathbb{P}\{X_{ij}(t)=k\}+o(1)\longrightarrow -\lambda t\sum_{k=1}^{r}p_k(1-z^k),\label{eq:sum-xij}
\end{align} where in the last step we used the condition \eqref{eq:TCCPP-triangular-3}. Hence using \eqref{eq:sum-xij}, from \eqref{eq:TCCPP-pgf-product} we have
\begin{align}
\log\mathbb{E}[z^{S_i(t)}]=\sum_{j=1}^{\infty}\log(1-x_{ij}(t))\longrightarrow -\lambda t \sum_{k=1}^{r}p_k(1-z^k)&=\lambda t \sum_{k=1}^{r}p_k(z^k-1).\nonumber
\end{align}
Therefore, 
\begin{equation}
    \mathbb{E}[z^{S_i(t)}]
    \longrightarrow
    \exp\{
        \lambda t
        \sum_{k=1}^{r}p_k(z^k-1)
    \}.\label{pgf_limit}
\end{equation}
Comparing \eqref{pgf_limit} with \eqref{LT_CPP}, we observe that \eqref{pgf_limit} is the pgf of a CPP with intensity $\lambda$ and $\mathbb{P}(X_j=k)=p_k$. Thus, $\{Y_i(t)\}_{t\ge 0}$ converges in distribution to a CPP. Now, consider the time-changed row sum process $\{Y_i(D_f(Y_g(t)))\}_{t\geq0}$.
For $-1\leq z\leq1$, its p.g.f. is given by
\begin{align}
\mathbb{E}[z^{Y_i(D_f(Y_g(t)))}]=\mathbb{E}[\mathbb{E}(z^{Y_i(D_f(Y_g(t)))}| D_f(Y_g(t)))].\label{pgf_time_changed}
\end{align}
Since $\{Y_i(t)\}_{t\geq0}$ is independent of
$\{D_f(Y_g(t))\}_{t\ge 0}$, conditioning on $D_f(Y_g(t))=s$, $Y_i(D_f(Y_g(t)))\overset{d}{=}Y_i(s)$.
Hence, using \eqref{pgf_limit},
\begin{align}
    \mathbb{E}[z^{Y_i(D_f(Y_g(t)))}]&\longrightarrow\mathbb{E}[\exp(\lambda D_f(Y_g(t))\sum_{k=1}^{r}p_k(z^k-1))]\label{pgf_time_changed_limit}.
\end{align}
Comparing \eqref{pgf_time_changed_limit} with the p.g.f. of the TCCPP
given in \eqref{p.g.f.TCCPP}, we conclude that
\begin{equation*}
    Y_i(D_f(Y_g(t)))
    \xrightarrow{d}
    \mathbb{M}_f^g(t),
    \quad i\to\infty.
\end{equation*}
\end{proof}

\begin{remark}\leavevmode
Theorem \ref{triangular_TCCPP} can be extended to the case where random variables $X_{ij}(t)$ may take values in an infinite state space provided $\sum_{k=1}^{\infty}p_k=1$ for the limiting TCCPP. For a triangular-array characterization of the compound Poisson distribution, see \citet{Zhang2016}. 
\end{remark}
\subsection{Characterization by waiting time distribution}

Next, we derive the waiting-time characterization of the TCCPP $\{\mathbb{M}_f^g(t)\}_{t\ge0}$. The inter-arrival times of a TCCPP have the same distribution as those of the underlying counting process $\{N_f^g(t):=N_f^g(t)\}_{t\ge 0}$. Being a time-changed Poisson process, the inter-arrival times of $\{N_f^g(t)\}_{t\ge 0}$ and hence those of the TCCPP are neither exponentially distributed nor independent. This implies that the distribution of their sum or the waiting time distribution of a TCCPP is difficult to obtain explicitly. 

Let $W_n^{f,g}$ be the waiting time until the $n$th arrival or the $n$th arrival time of the TCCPP. The distribution function of $W_n^{f,g}$ is given by 
\begin{equation}
    F_{W_n^{f,g}}(t)=\mathbb{P}(W_n^{f,g} \le t)=\mathbb{P}(\mathbb{M}_f^g(t) \ge n) =\sum_{i=n}^\infty p_f^g(i,t)=1 - \sum_{i=0}^{n-1} p_f^g(i,t).\label{eq:wt-cdf}
\end{equation} 
We first obtain a recursive relation satisfied by the arrival time densities of a TCCPP.
\begin{proposition}\label{proposition_waiting_time}
The $n$th arrival time of a TCCPP satisfies the following recursive relation
\begin{align*}
    \mathbb{P}\{W_n^{f,g}\in ds\}=\mathbb{P}\{W_{n-1}^{f,g}\in ds\}-ds\sum_{k=0}^{\infty}f_x^{*k}(n-1)\int_{0}^{\infty}\frac{e^{-\lambda x}(\lambda x)^{k}}{k!}\frac{d}{ds}\mathbb{P}\{D_f(Y_g(s))\in dx\}.
\end{align*}
\end{proposition}
\begin{proof} From \eqref{eq:wt-cdf} we have 
    \begin{align} 
   \mathbb{P}\{W_n^{f,g}\le s\}=\mathbb{P}\{\sum_{i=1}^{D_f(Y_g(s))}X_i\ge n\}
   &=\int_{0}^{\infty}\mathbb{P}\{\sum_{i=1}^{N_f^g(s)}X_i\ge n|D_f(Y_g(s))=x\}\mathbb{P}\{D_f(Y_g(s))\in dx\}\nonumber\\
   &=\int_{0}^{\infty}\mathbb{P}\{\sum_{i=1}^{N(x)}X_i\ge n\}\mathbb{P}\{D_f(Y_g(s))\in dx\}\nonumber\\
   &=\int_{0}^{\infty}\sum_{m=n}^{\infty}p(m,x)\mathbb{P}\{D_f(Y_g(s))\in dx\}.\nonumber
\end{align}
The density of the $n$th arrival time can be obtained as
\begin{align}
   \frac{\mathbb{P}\{W_n^{f,g}\in ds\}}{ds}=\frac{d}{ds}\mathbb{P}\{W_n^{f,g}\le s\}&=\frac{d}{ds}\int_{0}^{\infty}\sum_{m=n}^{\infty}\sum_{k=0}^{\infty}f_x^{*k}(m)\frac{e^{-\lambda x}(\lambda x)^{k}}{k!}\mathbb{P}\{D_f(Y_g(s))\in dx\}\quad(\text{using}~\eqref{pmf_CPP})\nonumber\\
    &=\frac{d}{ds}\int_{0}^{\infty}(1-\sum_{m=0}^{n-1}\sum_{k=0}^{\infty}f_x^{*k}(m)\frac{e^{-\lambda x}(\lambda x)^{k}}{k!})\mathbb{P}\{D_f(Y_g(s))\in dx\}\nonumber\\
    &=-\frac{d}{ds}\sum_{m=0}^{n-1}\sum_{k=0}^{\infty}f_x^{*k}(m)\int_{0}^{\infty}\frac{e^{-\lambda x}(\lambda x)^{k}}{k!}\mathbb{P}\{D_f(Y_g(s))\in dx\},\nonumber
\end{align}
which implies
\begin{align*}
    \mathbb{P}\{W_n^{f,g}\in ds\}=\mathbb{P}\{W_{n-1}^{f,g}\in ds\}-ds\sum_{k=0}^{\infty}f_x^{*k}(n-1)\int_{0}^{\infty}\frac{e^{-\lambda x}(\lambda x)^{k}}{k!}\frac{d}{ds}\mathbb{P}\{D_f(Y_g(s))\in dx\}.
\end{align*}
\end{proof}
The density of the first passage time or the time of first upcrossing of state $n$ of a TCCPP is the same as the density of its $n$th arrival time. 
Next, let $\mathbf{I}_n$ denote the $n\times n$ identity matrix and let $\mathbf{U}_n$ be the $n\times n$ shift matrix defined as
\begin{equation*}
    \mathbf{U}_n = \begin{pmatrix}
0 & 1 & 0 & \cdots & 0 \\
0 & 0 & 1 & \cdots & 0 \\
\vdots & \vdots & \ddots & \ddots & \vdots \\
0 & 0 & \cdots & 0 & 1 \\
0 & 0 & \cdots & 0 & 0
\end{pmatrix}.
\end{equation*}
Observe that $\mathbf{U}_n$ is a nilpotent matrix of degree $n$, that is, $\mathbf{U}_n^k = \mathbf{0}$ for all $k \ge n$. Furthermore, let $\boldsymbol{\beta} = (1, 0, \dots, 0)$ be a $1 \times n$ row vector, and let $\mathbf{1} = (1, 1, \dots, 1)^T$ be an $n \times 1$ column vector of ones. The following result provides a compact exponential matrix  representation for
 \eqref{eq:wt-cdf}.
\begin{theorem}\label{theorem_wating_time_distribution}
The distribution function of the $n$th arrival time of a TCCPP with intensity $\lambda$ and jump size distribution $\mathbb{P}(X_j=k)=p_k$ for $k\in\mathbb{N}$ is given by
\begin{equation*}
    F_{W_n^{f,g}}(t) = 1 - \boldsymbol{\beta} \, \mathbb{E}[ \exp(\mathbf{E}_n D_f(Y_g(t))) ] \mathbf{1},
\end{equation*}
where \begin{equation*}
    \mathbf{E}_n = \begin{pmatrix}
-\lambda & \lambda p_1 & \lambda p_2 & \cdots & \lambda p_{n-1}\\
0 & -\lambda & \lambda p_1& \cdots & \lambda p_{n-2} \\
0 & 0 & -\lambda & \cdots & \lambda p_{n-3}\\
\vdots & \vdots & \ddots & \ddots & \vdots \\
0 & 0 & \cdots & -\lambda & \lambda p_1\\
0 & 0 & \cdots & 0 & -\lambda
\end{pmatrix}. \qedhere
\end{equation*}
\end{theorem}

\begin{proof}
Since $\mathbf{E}_n$ is square and upper triangular with $-\lambda$ on its main diagonal, its only eigenvalue is $-\lambda$. Hence, its minimal polynomial is of the form $(x+\lambda)^k$ for some $k\leq n$, and therefore decomposes into linear factors. By Jordan Chevalley decomposition, $\mathbf{E}_n$ admits a unique decomposition into a diagonalizable part $\mathbf{D}$ and a nilpotent part $\mathbf{N}$ that commute as follows:
\begin{equation*}
    \mathbf{E}_n = -\lambda \mathbf{I}_n + \lambda \sum_{j=1}^{n-1} p_j \mathbf{U}_n^j,
\end{equation*}
where $\mathbf{D}=-\lambda\mathbf{I}_n$ and $\mathbf{N}=\lambda\sum_{j=1}^{n-1}p_j\mathbf{U}_n^j$.
Hence 
\begin{equation}
    \exp(\mathbf{E}_n t) = \exp(-\lambda t \mathbf{I}_n) \exp( \lambda t \sum_{j=1}^{n-1} p_j \mathbf{U}_n^j ) = e^{-\lambda t} \exp( \lambda t \sum_{j=1}^{n-1} p_j \mathbf{U}_n^j )=e^{-\lambda t} \sum_{m=0}^{\infty} \frac{(\lambda t)^m}{m!} ( \sum_{j=1}^{n-1} p_j \mathbf{U}_n^j )^m.\label{waiting time char eq1}
\end{equation}
Applying the multinomial theorem and the shift property $\mathbf{U}_n^a \mathbf{U}_n^b = \mathbf{U}_n^{a+b \pmod{n}}$, we have
\begin{equation}
    ( \sum_{j=1}^{n-1} p_j \mathbf{U}_n^j )^m = \sum_{\substack{c_1, \dots, c_{n-1} \ge 0 \\ \sum c_j = m}} \frac{m!}{c_1! \cdots c_{n-1}!} ( \prod_{j=1}^{n-1} p_j^{c_j} ) \mathbf{U}_n^{\sum_{j=1}^{n-1} j c_j \pmod{n}}.\label{multinomial_eq_1}
\end{equation}
Let $r \equiv \sum_{j=1}^{n-1} j c_j \pmod{n}$ so that $r\in\{0,1,\ldots,n-1\}$. Then \eqref{multinomial_eq_1} can be rewritten as
\begin{equation*}
    ( \sum_{j=1}^{n-1} p_j \mathbf{U}_n^j )^m = \sum_{r=0}^{n-1} ( \sum_{\substack{\sum c_j = m \\ \sum j c_j \equiv r \pmod{n}}} \frac{m!}{c_1! \cdots c_{n-1}!} \prod_{j=1}^{n-1} p_j^{c_j} ) \mathbf{U}_n^r.
\end{equation*}
The inner summation represents the total probability of a net displacement of $r$ after $m$ steps, which is the $m$-fold convolution of the jump distribution, denoted by $p_r^{*m}=f_X^{*m}(r)$. Thus, we obtain
\begin{equation}\label{waiting time char eq2}
    ( \sum_{j=1}^{n-1} p_j \mathbf{U}_n^j )^m = \sum_{r=0}^{n-1} p_r^{*m} \mathbf{U}_n^r.
\end{equation}
Substituting \eqref{waiting time char eq2} in \eqref{waiting time char eq1}, we have
\begin{equation*}
    \exp(\mathbf{E}_n t) = e^{-\lambda t}  \sum_{m=0}^{\infty} \frac{(\lambda t)^m}{m!} (\sum_{r=0}^{n-1} p_r^{*m}  \mathbf{U}_n^r)= \sum_{r=0}^{n-1} [ \sum_{m=0}^{\infty} e^{-\lambda t} \frac{(\lambda t)^m}{m!} f_X^{*m}(r) ] \mathbf{U}_n^r=\sum_{r=0}^{n-1} p(r, t) \mathbf{U}_n^r.
\end{equation*}
It follows that
\begin{align}
   \int_0^\infty \exp(\mathbf{E}_n x) h_f^g(x, t) \, dx&= \int_0^\infty ( \sum_{r=0}^{n-1} p(r, x) \mathbf{U}_n^r ) h_f^g(x, t) \, dx \nonumber \\
    \implies \mathbb{E}[ \exp(\mathbf{E}_n D_f(Y_g(t))) ] &= \sum_{r=0}^{n-1} ( \int_0^\infty p(r, x) h_f^g(x, t) \, dx ) \mathbf{U}_n^r \nonumber \\
    \implies \mathbb{E}[ \exp(\mathbf{E}_n D_f(Y_g(t))) ]&= \sum_{r=0}^{n-1} p_f^g(r, t) \mathbf{U}_n^r. \label{waiting time char eq4}
\end{align}
Explicitly, \eqref{waiting time char eq4} is the upper-triangular matrix
\begin{equation*}
     \mathbb{E}[ \exp(\mathbf{E}_n D_f(Y_g(t))) ] = \begin{pmatrix}
p_f^g(0,t) & p_f^g(1,t) & p_f^g(2,t) & \cdots & p_f^g(n-1,t) \\
0 & p_f^g(0,t) & p_f^g(1,t) & \cdots & p_f^g(n-2,t) \\
0 & 0 & p_f^g(0,t) & \cdots & p_f^g(n-3,t) \\
\vdots & \vdots & \ddots & \ddots & \vdots \\
0 & 0 & \cdots & p_f^g(0,t) & p_f^g(1,t) \\
0 & 0 & \cdots & 0 & p_f^g(0,t)
\end{pmatrix},
\end{equation*}
which implies
\begin{equation}\label{waiting time char eq5}
    \boldsymbol{\beta} \, \mathbb{E}[ \exp(\mathbf{E}_n D_f(Y_g(t))) ] \mathbf{1} = \sum_{r=0}^{n-1} p_f^g(r, t)=\mathbb{P}(\mathbb{M}_f^g(t) < n)=\mathbb{P}(W_n^{f,g} > t).
\end{equation}
Therefore, the distribution function of $W_n^{f,g}$ is given by
\begin{equation*}
    F_{W_n^{f,g}}(t) = 1 - \mathbb{P}(W_n^{f,g} > t) = 1 - \boldsymbol{\beta} \, \mathbb{E}[ \exp(\mathbf{E}_n D_f(Y_g(t))) ] \mathbf{1}. \qedhere
\end{equation*}
\end{proof}
\begin{corollary}\label{corollary_waiting_time}
For $\theta \in \mathbb{R}$, the characteristic function of the waiting time of a TCCPP is given by
\begin{equation*}
    \varphi_{W_n^{f,g}}(\theta) = \boldsymbol{\beta} [ \mathbf{I}_n + i\theta \int_0^\infty \exp(\mathbf{E}_n x) \hat{h}_f^g(x, \theta) \, dx ] \mathbf{1},
\end{equation*}
where $\hat{h}_f^g(x, \theta) = \int_0^\infty e^{i\theta t} h_f^g(x, t) \, dt$ denotes the Fourier transform of the density of $D_f(Y_g(t))$.
\end{corollary}
\begin{proof}
Since $f_{W_n^{f,g}}(t) = -\frac{\partial}{\partial t} \mathbb{P}(W_n^{f,g} > t)$, 
\begin{align}
    \varphi_{W_n^{f,g}}(\theta) = \int_0^\infty e^{i\theta t} f_{W_n^{f,g}}(t) \, dt&= [ -e^{i\theta t} \mathbb{P}(W_n^{f,g} > t) ]_0^\infty + i\theta \int_0^\infty e^{i\theta t} \mathbb{P}(W_n^{f,g} > t) \, dt\nonumber\\
    &=1 + i\theta \int_0^\infty e^{i\theta t} \mathbb{P}(W_n^{f,g} > t) \, dt\nonumber\\
    &= 1 + i\theta \int_0^\infty e^{i\theta t} \boldsymbol{\beta} [ \int_0^\infty \exp(\mathbf{E}_n x) h_f^g(x, t) \, dx ] \mathbf{1} \, dt.\quad (\text{using}~\eqref{waiting time char eq5})\nonumber
\end{align}
The boundary term in the integration-by-parts formula is equal to $1$ since $\mathbb{P}(W_n^{f,g}>t)\to0$ as $t\to\infty$, while $\mathbb{P}(W_n^{f,g}>0)=1$ due to the strictly positive jump sizes. Under the required integrability conditions, Fubini's theorem allows the order of integration to be interchanged. So 
\begin{align}
    \varphi_{W_n^{f,g}}(\theta) &= 1 + i\theta \boldsymbol{\beta} [ \int_0^\infty \exp(\mathbf{E}_n x) ( \int_0^\infty e^{i\theta t} h_f^g(x, t) \, dt ) dx ] \mathbf{1}=\boldsymbol{\beta} [ \mathbf{I}_n + i\theta \int_0^\infty \exp(\mathbf{E}_n x) \hat{h}_f^g(x, \theta) \, dx ] \mathbf{1}.\nonumber
\end{align}
\end{proof}
The following result provides the characterization of a TCCPP in terms of its $n$th arrival time.
\begin{theorem}\label{theorem_tail_probability}
A counting process is a TCCPP with rate $\lambda > 0$ and jump size distribution $\mathbb{P}(X_j=k)=p_k$ for $k\in\mathbb{N}$ if and only if the p.d.f. of its $n$th arrival time is given by  
\begin{equation}
    f_{W_n^{f,g}}(t)=-\frac{1}{2\pi} \int_0^{2\pi}\frac{e^{-i(n-1)\eta}}{1-e^{i\eta}} \frac{\partial}{\partial t}\mathbb{E}[ \exp( \lambda D_f(Y_g(t)) (\sum_{j=1}^\infty p_j e^{i\eta j} - 1) ) ] \, d\eta.\label{pdf_tail_distribution}
\end{equation}
\end{theorem}
\begin{proof}
Suppose $\{\mathbb{M}_f^g(t)\}_{t \ge 0}$ is a TCCPP. The tail probability of the TCCPP $\{\mathbb{M}_f^g(t)\}_{t\ge 0}$ is given by 
\begin{equation}
   \mathbb{P}(\mathbb{M}_f^g(t)\ge n)=\sum_{i=n}^\infty p_f^g(i,t)=\mathbb{P}(W_n^{f,g} \le t)=F_{W_n^{f,g}}(t).\nonumber
\end{equation}
For $|z| \le 1$, the generating function of the tail probability is 
    \begin{align}
        G_{W_n^{f,g}}(z,t) &= \sum_{n=1}^\infty F_{W_n^{f,g}}(t) z^n \label{generating_fun_definition}\\
    \implies (1-z)G_{W_n^{f,g}}(z,t) &= F_{W_1^{f,g}}(t)z + \sum_{n=2}^\infty (F_{W_n^{f,g}}(t) - F_{W_{n-1}^{f,g}}(t)) z^n \nonumber\\
    \implies (1-z)G_{W_n^{f,g}}(z,t)&= z(1 - p_f^g(0,t)) - \sum_{n=2}^\infty p_f^g(n-1,t) z^n \nonumber\\
    \implies (1-z)G_{W_n^{f,g}}(z,t)&= z - z \sum_{k=0}^\infty p_f^g(k,t) z^k \nonumber\\
    \implies G_{W_n^{f,g}}(z,t) &=\frac{z}{1-z}[1 - G_f^g(z,t)],\label{generating_fun_tail_prob}
    \end{align}
    where $G_f^g(z,t)$ is the p.g.f. of the TCCPP. 
    Differentiating both sides of \eqref{generating_fun_definition} with respect to $t$, we have
\begin{equation}\label{derivative_generating_fun}
    \frac{\partial}{\partial t}G_{W_n^{f,g}}(z,t) = \sum_{n=1}^\infty f_{W_n^{f,g}}(t) z^n.
\end{equation}
Now set $z=e^{i\eta}$ so $0\leq\eta\leq2\pi$. Applying the Fourier inversion formula to \eqref{derivative_generating_fun}, we obtain
\begin{align}
    f_{W_n^{f,g}}(t) &= \frac{1}{2\pi} \int_0^{2\pi} e^{-in\eta} \frac{\partial}{\partial t} G_{W_n^{f,g}}(e^{i\eta}, t) \, d\eta\nonumber\\
    &=\frac{1}{2\pi} \int_0^{2\pi} e^{-in\eta}(\frac{-e^{i\eta}}{1-e^{i\eta}} \frac{\partial}{\partial t} G_f^g(e^{i\eta},t)) \, d\eta\quad (\text{using}~\eqref{generating_fun_tail_prob})\nonumber\\
    &=-\frac{1}{2\pi} \int_0^{2\pi} \frac{e^{-i(n-1)\eta}}{1-e^{i\eta}} \frac{\partial}{\partial t}\mathbb{E}[ \exp( \lambda D_f(Y_g(t)) (\sum_{j=1}^\infty p_j e^{i\eta j} - 1) ) ] \, d\eta.\quad (\text{using}~\eqref{p.g.f.TCCPP})\nonumber
\end{align}
Conversely, assume that the $n$th arrival time of some counting process $\{A(t)\}_{t\ge 0}$ has density given by \eqref{pdf_tail_distribution}. 
The distribution function of the $n$th arrival time $W_n^A$ of $\{A(t)\}_{t\ge 0}$ is therefore given by
\begin{align}
\mathbb P(W_n^{A}\le t)
&=-\frac{1}{2\pi}\int_0^{2\pi}
\frac{e^{-i(n-1)\eta}}{1-e^{i\eta}}[\mathbb E[
\exp\{
\lambda D_f(Y_g(t))
(\sum_{j=1}^{\infty}p_je^{i\eta j}-1)
\}
]-1]d\eta.\nonumber
\end{align}
Since $\{W_n^A\le t\}=\{A(t)\ge n\}$, the generating function of the tail probability of $\{A(t)\}_{t\ge 0}$ for $|z|\le 1$ is
\begin{align}
Q(z,t)=\sum_{n=1}^{\infty}
\mathbb P(A(t)\ge n)z^n
&=-\frac{1}{2\pi}\int_0^{2\pi}
\frac{1}{1-e^{i\eta}}[\mathbb E[
\exp\{
\lambda D_f(Y_g(t))
(\sum_{j=1}^{\infty}p_je^{i\eta j}-1)
\}
]-1](\sum_{n=1}^{\infty}z^n e^{-i(n-1)\eta})d\eta\nonumber\\
&=-\frac{1}{2\pi}\int_0^{2\pi}
\frac{1}{1-e^{i\eta}}[\mathbb E[
\exp\{
\lambda D_f(Y_g(t))
(\sum_{j=1}^{\infty}p_je^{i\eta j}-1)
\}
]-1]\frac{z}{1-ze^{-i\eta}}d\eta.\nonumber
\end{align}
Substituting $w=e^{i\eta}$ gives $dw=iw\,d\eta$, so that
$d\eta=\frac{dw}{iw}$. Thus, the integral over $\eta\in[0,2\pi]$
becomes a contour integral along the unit circle $|w|=1$. It follows that
\begin{align}
    Q(z,t)&=-\frac{1}{2\pi}\oint_{|w|=1}
\frac{1}{1-w}[\mathbb E[
\exp\{
\lambda D_f(Y_g(t))
(\sum_{j=1}^{\infty}p_jw^{j}-1)
\}
]-1]\frac{z}{1-z/w}\frac{dw}{iw}\nonumber\\
&=-\frac{z}{2\pi i}\oint_{|w|=1}
\frac{1}{(1-w)(w-z)}[\mathbb E[
\exp\{
\lambda D_f(Y_g(t))
(\sum_{j=1}^{\infty}p_jw^{j}-1)
\}
]-1]dw\nonumber\\
&=-\frac{z}{1-z}\frac{1}{2\pi i}\oint_{|w|=1}
 ( \frac{1}{1-w} + \frac{1}{w-z} )[\mathbb E[
\exp\{
\lambda D_f(Y_g(t))
(\sum_{j=1}^{\infty}p_jw^{j}-1)
\}
]-1]dw.\label{split_integrals}
\end{align}
Splitting the integral into two parts, we evaluate each term separately. Since $\sum_{j=1}^{\infty}p_j=1$, the bracketed term vanishes at $w=1$, making the singularity in the first integral removable so that it evaluates to zero by Cauchy's integral theorem. For the second integral, the simple pole at $w=z$ lies inside the unit circle for $|z|<1$. Hence applying Cauchy's integral formula, \eqref{split_integrals} reduces to
\begin{align}
    Q(z,t)&=\frac{z}{z-1}(\mathbb{E}[
\exp\{
\lambda D_f(Y_g(t))
(\sum_{j=1}^{\infty}p_jz^{j}-1)
\}
]-1).\label{waiting_time_generating_function_A}
\end{align}
Now, the p.g.f. of $\{A(t)\}_{t\ge 0}$ is given by
\begin{align}
    G_A(z,t)&=\mathbb{P}(A(t)=0)
+\sum_{n=1}^{\infty}\mathbb{P}(A(t)=n)z^n\nonumber\\
&=1-\mathbb{P}(A(t)\geq1)
+\sum_{n=1}^{\infty}\mathbb{P}(A(t)\geq n)z^n-\frac{1}{z}[
\sum_{n=1}^{\infty}\mathbb{P}(A(t)\geq n)z^n
-\mathbb{P}(A(t)\geq1)z
]\nonumber\\
&=1+\frac{z-1}{z}Q(z,t)\nonumber\\
&= \mathbb{E}[ \exp\{ \lambda D_f(Y_g(t)) (\sum_{j=1}^{\infty}p_jz^j-1) \} ],\quad(\text{using}~\eqref{waiting_time_generating_function_A})\nonumber
\end{align}
which is the p.g.f. of a TCCPP with rate $\lambda$ and jump size probabilities $p_j$ for $j\in\mathbb{N}$. Hence, $\{A(t)\}_{t\ge 0}$ is a TCCPP.
\end{proof}
\begin{remark}
If $D_f(Y_g(t))=t$ a.s. (TCCPP reduces to CPP), then Theorem \ref{theorem_wating_time_distribution} and Theorem \ref{theorem_tail_probability} reduce to Theorem 6.1 and Theorem 6.2 respectively of \citet{Zhang2016}.    
\end{remark}
Finally, if the TCCPP is defined as in Definition~\ref{def:TCCPP-2}, then all results established for the process in Definition~\ref{def:TCCPP-1} in this paper remain valid when \( D_f(Y_g(t)) \) is replaced by \( Y_g(D_f(t)) \). The only exception is the law of the iterated logarithm stated in Theorem~\ref{lil_TCCPP}, which does not extend in the same form to this setting. This is due to the fact that Lemma~\ref{lil.Subordinator} used in its proof, applies specifically to subordinators.

\section*{Declarations}
\subsection*{Conflict of Interest} The authors declare that they have no conflict of interest.
\subsection*{Data availability} No data was used for the research described in this article.
\subsection*{Ethics Approval} Not applicable for the research described in this article. 
\subsection*{Funding} The authors received no specific funding for this research.
\subsection*{Declaration of generative AI and AI-assisted technologies in the manuscript preparation process} The authors declare that no generative AI or AI-assisted technologies were used in preparing this article.
\bibliographystyle{apalike}
\bibliography{TCCPP}

\end{document}